\documentclass[reqno]{amsart}
\usepackage{hyperref}
\usepackage{fancyvrb}
\usepackage{amsmath}
\usepackage{amsfonts}
\usepackage{amsthm}
\usepackage{amssymb}
\usepackage{array}
\usepackage{cases}
\usepackage[hypcap=false]{caption}
\usepackage{xcolor}
\hypersetup{colorlinks=true,linkcolor=blue!60!black,citecolor=green!40!black,urlcolor=blue!60!black}
\usepackage{tikz}
\usetikzlibrary{calc}
\usepackage{centernot}
\usepackage{mathtools}
\usepackage{graphicx}
\usepackage[margin=1in]{geometry}
\usepackage[makeroom]{cancel}
\usepackage{stmaryrd}
\usepackage{stackengine}

\newtheorem{lemma}{Lemma}[subsection]
\newtheorem{corollary}[lemma]{Corollary}

\newtheorem{proposition}[lemma]{Proposition}
\newtheorem{theorem}{Theorem}[section]

\theoremstyle{definition}
\newtheorem{definition}[lemma]{Definition}
\newtheorem{example}[lemma]{Example}

\theoremstyle{remark}

\newtheorem*{remark}{Remark}

\numberwithin{equation}{subsection}

\DeclareMathOperator{\sch}{\mathfrak{S}}
\newcommand{\tom}[1]{\xrightarrow{#1}}
\newcommand{\Tom}[1]{\xRightarrow{#1}}

\newcommand{\cperm}[1]{\llbracket #1\rrbracket}
\newcommand{\cpcf}[5]{d^{#1,#2}_{#3}\left(#4\mid #5\right)}
\newcommand{\smpfu}{\mathsf{\Pi}}
\newcommand{\smpr}[3]{\ensuremath{}\smpfu\left( #1 \mid #2_{[#3]}\right)}
\newcommand{\smpre}[2]{\ensuremath{}\smpfu\left( #1 \mid #2\right)}
\newcommand{\downvar}[1]{\stackrel{#1}{\searrow}}

\newcommand{\shifthom}{\tau}

\newcommand{\shiftvby}[1]{\shifthom^{#1}}

\makeatletter
\newlength{\negph@wd}
\DeclareRobustCommand{\negphantom}[1]{%
	\ifmmode
	\mathpalette\negph@math{#1}%
	\else
	\negph@do{#1}%
	\fi
}
\newcommand{\negph@math}[2]{\negph@do{$\m@th#1#2$}}
\newcommand{\negph@do}[1]{%
	\settowidth{\negph@wd}{#1}%
	\hspace*{-\negph@wd}%
}
\makeatother

\makeatletter
\newlength{\negvph@wd}
\DeclareRobustCommand{\negvphantom}[1]{%
	\ifmmode
	\mathpalette\negvph@math{#1}%
	\else
	\negvph@do{#1}%
	\fi
}
\newcommand{\negvph@math}[2]{\negvph@do{$\m@th#1#2$}}
\newcommand{\negvph@do}[1]{%
	\settowidth{\negvph@wd}{#1}%
	\hspace*{-\negvph@wd}%
}
\makeatother
\newcommand*\circled[1]{%
	\tikz[baseline=(char.base)]{
		\node[shape=circle,draw,inner sep=0.1pt] (char) {#1};
	}%
}
\newcommand\mycircled[1]{\makebox[0pt]{\circled{#1}}}
\newcommand{\tb}{\textbullet}
\newcommand{\xr}[1]{\ensuremath{{}\xrightarrow{[#1]}{}}}
\newcommand{\rd}[1]{\ensuremath{\mathcolor{red}{\mathbf{#1}}}}
\newcommand{\mc}[1]{\ensuremath{\mycircled{#1}}}
\newcommand{\desc}{\mathrm{Desc}}%

\title[Dominant permutations and Schubert polynomials]{Dominant permutations, positive mixed-variable formulas, and product--splitting equivalence for (double) Schubert polynomials}
\author{Matthew J. Samuel}
\email{matthew.samuel.math@gmail.com}

\newcommand{\dom}{\mathfrak{d}}
\newcommand{\ldom}{\dom^*}
\newcommand{\code}{\mathfrak{c}}

\begin{document}
\newlength{\heightoff}
\settoheight{\heightoff}{$f$}

\begin{abstract}
	We present new positive combinatorial formulas for triple Schubert coefficients in mixed-variable products of double Schubert polynomials, $\sch_u(x;y)\sch_v(x;z)$, expressed in the $\sch_w(x;y)$ basis, in the case where $u$ is dominant. This result complements earlier results where $v$ was assumed dominant. To generalize further, we introduce a relation on permutations that we call ``domination,'' under which dominant permutations dominate all others. When $u$ dominates $w$, we obtain an explicit positive formula for $c^w_{u,v}(y;z)$, which in particular yields positivity when $u$ and $v$ have separated descents.

	In addition to these new positivity results, we establish a product--splitting equivalence for double Schubert polynomials, relating triple Schubert coefficients to splitting coefficients. This equivalence highlights structural symmetries and provides reduction methods for computing coefficients that are otherwise difficult to access. Taken together, these results extend the range of double Schubert polynomial products for which positivity can be established combinatorially and yield explicit formulas that can be implemented computationally.

\end{abstract}

\maketitle

\newcommand{\dsch}{\ensuremath{\Xi}}
\newcommand{\coma}{\mathcal{A}}
\newcommand{\dcoma}{\mathcal{D}}

\tableofcontents

\section*{Introduction}

Double Schubert polynomials $\sch_u(x;y)$, indexed by permutations, constitute a central basis in Schubert calculus. Their multiplication with respect to different sets of coefficient variables defines triple Schubert structure constants, also called triple Schubert coefficients, via the relation
\begin{equation*}
	\sch_u(x;y)\,\sch_v(x;z)\;=\;\sum_w c^w_{u,v}(y;z)\,\sch_w(x;y).
\end{equation*}
The central objective in this area is to obtain explicit combinatorial formulas that exhibit positivity, in analogy with the Littlewood--Richardson rule for Schur polynomials. Since the triple structure constants are polynomials, their positivity is understood in terms of suitable differences of coefficient variables, rather than ordinary monomials in $y$ and $z$.

Gao and Xiong proved the positivity conjecture of \cite[Conjecture~1.1]{samuelmolev} in their FPSAC 2026 paper \emph{Positivity in triple Schubert calculus} \cite[Theorem~1.1]{gao2025grahampositivitytripleschubert}. With our ordering of the factors and coefficient variables, their theorem states that every triple Schubert coefficient $c_{u,v}^w(y;z)$ is a polynomial in the differences $y_i-z_j$ with nonnegative integer coefficients. Their geometric proof uses a refinement of Graham's positivity theorem \cite{graham2001positivity}, but does not provide a general positive combinatorial formula. Positivity in arbitrary differences $y_i-z_j$ need not make Graham positivity in the simple roots $y_{i+1}-y_i$ manifest after specializing $z=y$.

In this paper, we give explicit positive combinatorial formulas in broad new cases, complementing Gao and Xiong's geometric positivity theorem. We introduce a relation on permutations that we call ``domination,'' extending the classical notion of dominant permutation. We prove that if $u$ dominates $w$, then the coefficient $c^w_{u,v}(y;z)$ has a positive combinatorial formula. This complements the earlier formulas of \cite{samuelmolev} for dominant $v$ and yields positive combinatorial formulas when $u$ and $v$ have separated descents in either orientation.

We also establish an equivalence between products and splitting coefficients of double Schubert polynomials. This equivalence enables us to compute products from splitting coefficients and conversely, thereby exposing structural symmetries and providing practical methods to reduce complex cases to simpler ones. Bergeron and Sottile identified the splitting coefficients of a Schubert polynomial $\sch_w$ with structure constants of products in which $\sch_w$ is one of the factors \cite[Theorem~4.6.1]{coprod}. Our equivalence identifies them with structure constants of different products, whose two factors correspond to the two Schubert polynomials in the splitting and whose product term is determined by $w$; it is new already for ordinary Schubert polynomials (see the remark following Theorem \ref{theorem:doublecoprodpos}). The equivalence also gives a framework for organizing the difficulty of multiplication problems in terms of splitting-coefficient computations.

\subsection*{Software used}

All computations in this work were carried out using \texttt{schubmult} \cite{schubmult}, a free Python package developed by the author. The package is available on PyPI and comes with command-line executables. It supports:
\begin{itemize}
    \item Products of ordinary Schubert polynomials
    \item Products of double Schubert polynomials (both same-variable and mixed-variable)
    \item Splitting coefficients of Schubert polynomials
    \item Products of quantum and quantum double Schubert polynomials (including same-variable, mixed-variable, and parabolic cases)
    \item Conjectural analogues for Grothendieck, double Grothendieck, and quantum double Grothendieck polynomials
\end{itemize}
Questions, bug reports, or other inquiries about the package can be directed to \href{mailto:schubmult@gmail.com}{schubmult@gmail.com}.

\section{The domination formula} \label{section:dualpieri}

\subsection{Statement of the formula with examples}

\begin{definition}[$S_\infty$ as a Coxeter group, window notation]
	For standard results and notation regarding Coxeter groups, we refer the reader to \cite{combcox} or \cite{humphreys1992reflection}. $S_\infty$ will refer to the infinite symmetric group of permutations $\mathbb N\to\mathbb N$ fixing all but finitely many elements, and $S_n$ is the subgroup of all $u\in S_\infty$ such that $u(i)=i$ for all $i>n$. We will use window notation, where we express a permutation $u$ as a permutation of the first $n$ positive integers (enclosed by brackets $[]$) such that if we write
	$$u=[a_1,a_2,\ldots,a_n]$$
	where $a_i\neq a_j$ if $i\neq j$ and $1\leq a_i\leq n$ for all $i$, then $u(i)=a_i$ for $1\leq i\leq n$, and $u(i)=i$ if $i>n$. An alternative notation, code notation $\cperm{\cdots}$, will be defined in Definition \ref{definition:codestuff}.

	We use the standard notation of $s_i$ as the adjacent transposition such that right multiplication exchanges the elements at indices $i$ and $i+1$, and $t_{ij}$ is the transposition such that right multiplication exchanges the elements at indices $i$ and $j$. These are the simple reflections (resp. reflections) in a Coxeter presentation of $S_\infty$.
\end{definition}

\begin{definition}[Separated descents] \label{definition:sepdesc}
	Two permutations $u$ and $v$ are said to have \emph{separated descents} if there exists a positive integer $p$ such that $\ell(us_i) > \ell(u)$ for all $i < p$ and $\ell(vs_j) > \ell(v)$ for all $j>p$. It is important to note that this definition is not symmetric in $u$ and $v$. The formulas in \cite{samuelmolev} and \cite{fan2025bumpless} establish positivity of $c_{u,v}^w(y;z)$ in the separated descents case as defined above. A positive combinatorial formula for $c_{v,u}^w(y;z)$ is provided in this article.
\end{definition}

\begin{definition}[Special group elements {$c[p,q]$}, {$d[p,q]$}, {$r[p,q]$}, the homomorphism $\sigma_p$] \label{definition:specelem}
	For $1\leq p\leq k$, define the cycle
	$$c[p,k]=s_{k+1-p}s_{k+2-p}\cdots s_k$$
	For $p,k\geq 1$, also define
	$$r[p,k]=s_{k+p-1}s_{k+p-2}\cdots s_k$$
	We use the same notation for these as in \cite{sottile}, and $c[p,k]$ corresponds to elementary symmetric polynomials, while $r[p,k]$ corresponds to complete symmetric polynomials. The letters $c$ and $r$ stand for ``column'' and ``row,'' respectively, referring to the shape of the indexing partition of the corresponding Schur polynomial.

	For $q\geq 0$ define
	$$d[p,q]=s_ps_{p+1}\cdots s_{p+q-1}$$
	(where by convention $d[p,0]=1$) and define
	$$d[q]=d[1,q]$$
	For $q>0$, we have $d[p,q]=c[q,p+q-1]$, so that this is merely a convenient way to re-express the $c$ permutations. If $0\leq q_1\leq q_2$, then $d[p,q_1]^{-1}d[p,q_2]=d[p+q_1,q_2-q_1]$.

		Define a partial homomorphism $\sigma_p$ for a not necessarily positive integer $p$ by $\sigma_p(s_i)=s_{i+p}$ whenever $i+p>0$. This is defined for all $w$ such that $w=s_{i_1}\cdots s_{i_m}$ with $i_j+p>0$ for all $j$. $\sigma_1$ will be used often, so $\sigma$ with no subscript will indicate $\sigma_1$ to avoid clutter.
	\end{definition}

	\begin{definition}[Sottile's Pieri relations $\tom{k},\Tom{k}$] \label{definition:pierisymb}
		Following \cite{sottile}, with slightly different notation, for permutations $u,w$ we declare that $u\tom{k} w$ if there exist transpositions $t_{a_1b_1},\ldots, t_{a_pb_p}$ for some $p$ with $0\leq p\leq k$ such that:
		\begin{enumerate}
			\item \label{item:it1} $a_i\leq k<b_i$ for all $i$.
			\item \label{item:it2} $a_i\neq a_j$ if $i\neq j$.
			\item \label{item:it3} $\ell(ut_{a_1b_1}\cdots t_{a_ib_i})=\ell(u)+i$ for all $1\leq i\leq p$.
			\item \label{item:it4} $w=ut_{a_1b_1}\cdots t_{a_pb_p}$.
		\end{enumerate}
		This definition is used in Sottile's Pieri formula in the reference and is central to the development of our formula. The second relation, also due to Sottile, is defined by declaring that $u\Tom{k} w$ if there exist transpositions $t_{a_1b_1},\ldots,t_{a_pb_p}$ for some $p\geq 0$ such that (\ref{item:it1}), (\ref{item:it3}), and (\ref{item:it4}) hold, but instead of (\ref{item:it2}) we have that $b_i\neq b_j$ if $i\neq j$. Note that, unlike for $\tom{k}$, $p$ can be greater than $k$ with the relation still holding.

	\end{definition}

	We next introduce the code of a permutation, a compact representation that we will use together with the Pieri relations.

	\begin{definition}[$\mathfrak{c}(w)$, $\mathfrak{c}^*(w)$, dominant permutations, the dominant approximation $\dom(v)$, code representation $\cperm{\cdots}$] \label{definition:codestuff}

        Given $w\in S_\infty$, we define the \emph{code} $\code(w)$ (also known as the Lehmer code) of $w$ as follows:
		$$\code_i(w)=\#\{j> i\mid w(i)>w(j)\}$$
		In other words, $\code_i(w)$ is the number of inversions $(i,j)$ of $w$, with $i$ fixed. We also define the dual code
		$$\code^*(w)=\code(w^{-1})$$
		which is simply the code of the inverse. Permutation codes already appear, for instance, in \cite{laisant1888numeration}.

		Some useful facts about the code, which we will leave to the reader, are as follows: $u(i)<u(i+1)$ if and only if $\code_i(u)\leq\code_{i+1}(u)$, and in that case
		$$\code_j(us_i)=\begin{cases}
			\code_j(u)&\mbox{ if }j\notin \{i,i+1\}\\
			\code_{i+1}(u)+1&\mbox{ if }j=i\\
			\code_i(u)&\mbox{ if }j=i+1
		\end{cases}$$
		Furthermore, $\code^*_i(w)$ is equal to the number of inversions $(i,j)$ of $w^{-1}$, which is to say that $w^{-1}(i)>w^{-1}(j)$. This means that, in terms of $w$, if $q$ is the index such that $w(q)=i$, then $\code^*_i(w)$ is the number of indices $p<q$ such that $w(p)>i$. It follows easily that
		$$w(\code^*_1(w)+1)=1$$
		The inverse function of the code will be denoted by $\cperm{\cdots}$, meaning $\cperm{\code(w)}=w$. For example, since
		$$\code([1,3,5,2,4,7,6]) = (0,1,2,0,0,1)$$
		we have
		$$\cperm{0,1,2,0,0,1}=[1,3,5,2,4,7,6]$$

		A permutation $\mu$ is said to be a \emph{dominant permutation} if $\code(\mu)$ is a partition. Equivalently, $\mu$ avoids the pattern $132$: there are no indices $i<j<k$ such that $\mu(i)<\mu(k)<\mu(j)$ \cite{reif}.

		For a permutation $v$, we define a dominant permutation $\dom(v)$ following \cite{samuelmolev}, recursively as follows. If $v$ is dominant, then $\dom(v)=v$. Otherwise, let $i$ be the maximal index such that $\code_i(v)<\code_{i+1}(v)$; then set
		$$\dom(v)=\dom(vs_i)$$
		If $v\in S_n$, each step stays in $S_n$ and increases the length by $1$, so the process terminates. The prescribed choice of index makes $\dom(v)$ well defined.
	\end{definition}

	\begin{definition}[$\varphi_{i,n}$, $\mathcal{D}_i(v)$, $Q_i(v',v)$] \label{definition:polysequence}
		For $v'\in S_n$ and $1\leq i\leq n+1$, define
		$$\varphi_{i,n}(v')(j)=\begin{cases}
			v'(j)&\mbox{ if }j<i\\
			n+1&\mbox{ if }j=i\\
			v'(j-1)&\mbox{ if }i<j\leq n+1\\
			j&\mbox{ if }j>n+1
		\end{cases}$$
		For $v,v'\in S_\infty$ and an integer $i\geq 1$, we define a relation $\downvar{i}$ by declaring that $v\downvar{i} v'$ if, for all sufficiently large $n$ with $v,v'\in S_n$ and $n\geq i$, we have
		$$v\Tom{i}\varphi_{i,n}(v')$$
		Thus $v$ satisfies the relation $\Tom{i}$ with respect to the permutation obtained from $v'$ by inserting $n+1$ at position $i$, after padding the permutations with fixed points as necessary. This concept was introduced by Bergeron and Sottile in \cite{bsskew}.

		If $v\downvar{i} v'$, we define a set of integers $Q_i(v',v)$ by
		$$Q_i(v',v)=\{v(j)\mid j>i\mbox{ and }v'(j-1)=v(j)\}$$
		Given the fact that any element of $S_\infty$ fixes all but finitely many positive integers, it follows that $Q_i(v',v)$ is a finite set.

		We then define $\mathcal{D}_i(v)$ to be the set of all permutations $v'\in S_\infty$ such that $v\downvar{i} v'$. The permutations in $\mathcal{D}_i(v)$ arise when pulling a variable out of a Schubert polynomial and expressing the coefficients of the powers of this variable as Schubert polynomials in the remaining variables, as is done in \cite{coprod}, from which this definition essentially comes. The associated set $Q_i(v',v)$ is specific to the case of double Schubert polynomials.
	\end{definition}

	\begin{example}
		The permutation $\mu=[5,4,6,3,1,2]$ has code $(4,3,3,2)$ and hence is dominant. We then have that $\code^*(\mu)=(4,4,3,1)$.

		Let $v=[4,3,1,5,2]$. Then $\code(v)=(3,2,0,1)$, hence $\dom(v)=vs_3=[4,3,5,1,2]$. We have that $\code(\dom(v))=(3,2,2)$ and $\code^*(\dom(v))=(3,3,1)$. We have that $v^{-1}=[3,5,2,1,4]$, $\code(v^{-1})=(2,3,1)$, hence $\dom(v^{-1})=v^{-1}s_1=[5,3,2,1,4]$. $\code(\dom(v^{-1}))=(4,2,1)$ and $\code^*(\dom(v^{-1}))=(3,2,1,1)$. We have that $v\downvar{3}[4,3,1,2]$ and $v\downvar{3}[4,3,2,1]$, and this list is exhaustive.
		We then have that
		$$Q_3([4,3,1,2],v)=\{2\}$$
		and $Q_3([4,3,2,1],v)=\emptyset$.
	\end{example}

	\begin{definition}[$\phi$, $\phi_m$, \emph{$k$-dominates}, \emph{dominates}] \label{definition:domination}
		A permutation $u$ is said to \emph{$1$-dominate} the permutation $w$ if for each index $j$ with $\code_1^*(u)<j\leq \code_1^*(w)$ we have $u(j)<u(j+1)$. Equivalently, $u$ has no descents strictly to the right of position $\code_1^*(u)$ up to position $\code_1^*(w)$.

		For $i\geq 0$, we define a function $\phi_i:S_\infty\to S_\infty$ by declaring that
		$$\code^*_j(\phi_i(u))=\code^*_{j+i}(u)$$
		for all $j$, essentially ``chopping off'' the first $i$ elements of the dual code. We recursively say that $u$ $k$-\emph{dominates} $w$ for $k>1$ if $u$ $(k-1)$-dominates $w$ and $\phi_{k-1}(u)$ $1$-dominates $\phi_{k-1}(w)$. We then say that $u$ \emph{dominates} $w$ if $u$ $k$-dominates $w$ for all $k\geq 1$. For simplification of notation, we write
		$$\phi=\phi_1$$
	\end{definition}

	\begin{example}
		The permutation $u=[2,3,5,1,4]$ dominates the permutation $w=[2,6,3,4,7,1,5]$ but is not dominant. To see that $u$ dominates $w$, note that $\code_1^*(u)=3$, $\code_1^*(w)=5$, and $u(4)<u(5)<u(6)$, so $u$ $1$-dominates $w$. We have that $\phi(u)=[1,2,4,3]$ and $\phi(w)=[1,5,2,3,6,4]$. Clearly $\phi(u)$ $1$-dominates $\phi(w)$ since $\code_1^*(\phi(u))=\code_1^*(\phi(w))$. We have that $\phi_2(u)=[1,3,2]$ and $\phi_2(w)=[4,1,2,5,3]$. Since $\code_1^*(\phi_2(u))=0$, $\code_1^*(\phi_2(w))=1$, and $\phi_2(u)(1)<\phi_2(u)(2)$, we have that $u$ $3$-dominates $w$. It is indeed true that $u$ dominates $w$.
	\end{example}

	\begin{definition}[Double Schubert polynomials $\sch_u(x;y)$, divided difference operators $\partial^v$, $\partial_v^w$, $c_{u,v}^w(y;z)$] \label{definition:schubert}
		We denote the double Schubert polynomial indexed by $u\in S_\infty$ in the variables $x$ and $y$ by $\sch_u(x;y)$. We define coefficients $c_{u,v}^w(y;z)$ by
		$$\sch_u(x;y)\sch_v(x;z)=\sum_{w\in S_\infty} c_{u,v}^w(y;z)\sch_w(x;y)$$
		Divided difference operators are indexed by elements of $S_\infty$ written as superscripts, for example $\partial^u$. We define skew divided difference operators $\partial_v^w$ by requiring that, for polynomials $P,Q$,
		$$\partial^w(PQ)=\sum_{v\in S_\infty}\partial^v(P)\partial_v^w(Q)$$
		This differs from the definition of skew divided difference operators in \cite{notes}, but agrees with that in \cite{samuelleibniz}, and in our opinion is more natural. We have that
		$$c_{u,v}^w(x;z)=\partial_u^w(\sch_v(x;z))$$
		An explanation of this fact can be found in \cite{samuelmolev}.
	\end{definition}

	\begin{definition}[Polynomial variable omission $x^{(i)}$, subscript substitution $y_A$, index shift function $\shifthom$] \label{definition:pv}
		Define
		$$x^{(i)}=(x_1,\ldots,\widehat{x_i},\ldots)$$
		For a sequence of multiple indices in the superscript, we omit all of them.

		For a set of positive integers $A$ and a sequence of indeterminates such as $y = (y_i)_{i=1}^{\infty}$, we denote
		$$y_A = \{y_a\mid a\in A\}$$
		While this is expressed as a set, the implication is that the possibly nonconsecutive indices occur in increasing order.

		We define a homomorphism $\shifthom:\mathbb{Z}[x,y,z]\to\mathbb{Z}[x,y,z]$ by specializing $\shifthom(x_j)=x_{j+1}$, $\shifthom(y_j)=y_j$ and $\shifthom(z_j)=z_j$. We denote the composition of $\shifthom$ with itself $i$ times by $\shiftvby{i}$. We also apply this notation to sequences: for a sequence $a=(a_1,a_2,\ldots)$ we write $\shiftvby{i}a=(a_{i+1},a_{i+2},\ldots)$. For example, $\sch_v(\shiftvby{i}y;z)$ is obtained from $\sch_v(x;z)$ by substituting $y_{j+i}$ for $x_j$ for all $j$.

	\end{definition}

	\begin{definition}[$H_p(x;y)$, $E_p(x;y)$, $\smpre{a}{B}$] \label{definition:hp}
		Define
		$$H_p(x;y)=\prod_{i=1}^p (x_1-y_i)$$
		This is actually a double Schubert polynomial, namely
		$$H_p(x;y)=\sch_{r[p,1]}(x;y)$$
		Define also
		$$E_p(x;y)=\prod_{i=1}^p (x_i-y_1)$$

		Let $\mathcal{F}$ denote the set of finite subsets of $\mathbb{Z}[x,y,z]$. We define a function
		$$\smpfu:\mathbb{Z}[x,y,z]\times \mathcal{F}\to\mathbb{Z}[x,y,z]$$
		that is $H_p(x;y)$ thinly disguised by substitution, by
		$$\smpre{a}{B} = \prod_{b\in B}(a-b)$$
		Typically, $a$ will be a single indeterminate and $B$ will be a set or sequence of indeterminates, rather than arbitrary elements of the ring. For example,
		$$\smpre{x_1}{y_{[p]}} = \prod_{i=1}^p(x_1-y_i) = H_p(x;y)$$
		$$\smpre{-y_1}{(-x)_{[p]}} = \prod_{i=1}^p(x_i-y_1) = E_p(x;y)$$
		By a mild abuse of notation, we may write such expressions as
		$$\smpre{x_i}{y_1,y_3,y_4} = (x_i-y_1)(x_i-y_3)(x_i-y_4)$$

	\end{definition}

 We may now state our domination formula. The notation is given in Definitions \ref{definition:specelem}--\ref{definition:hp}.

\begin{theorem}[The domination formula] \label{theorem:dualpieri}
Let $u,v,w\in S_\infty$ and suppose $u$ dominates $w$. Then $c_{u,v}^w(y;z)=0$ unless $\code_j^*(w)\geq \code_j^*(u)$ for all $j$. If that condition is satisfied, let $m$ be the maximal index such that $\code_m^*(u)\neq 0$, with $m=0$ when $u=1$. For each $1\leq k\leq m$ define
$$d_k=d[\code_k^*(u)+1,\code^*_k(w)-\code_k^*(u)]$$
Then
$$c_{u,v}^w(y;z)=\sum_{v_1,\ldots,v_m} \left(\prod_{i=1}^m\smpre{y_i}{z_{A_i}}\right)\sch_{v_m\phi_m(w)^{-1}}(\shiftvby{m}y;z)$$
where $v_0=v$, $\ell(v_m\phi_m(w)^{-1})=\ell(v_m)-\ell(\phi_m(w))$, and $v_k$ ranges over all elements such that
$$v_{k-1}d_k^{-1}\downvar{\code_k^*(u)+1}v_k$$
and
$$\ell(v_{k-1}d_k^{-1})=\ell(v_{k-1})-\code_k^*(w)+\code_k^*(u)$$
and we define
$$A_k = Q_{\code_k^*(u)+1}(v_k,v_{k-1}d_k^{-1})$$
for all $k$.

In particular, this holds for all $v$ and $w$ when $u$ is dominant (see Lemma \ref{lemma:domdom}).
\end{theorem}

We illustrate the terms in the sum with diagrams, which are conceptually similar to the diagrams in \cite{samuelmolev}. In each diagram, the relevant permutations are arranged column-wise. Circled entries represent elements of the $Q$ set (Definition \ref{definition:polysequence}). Visually, these correspond to numbers that do not change from column to column.

 More specifically, each diagram is a rectangle containing the permutations $v_i$, $v_id_{i+1}^{-1}$, and $v_m\phi_m(w)^{-1}$ in sequence. If $v_i\neq v_id_{i+1}^{-1}$, then there will be an arrow from $v_i$ to $v_id_{i+1}^{-1}$ labeled with $d_{i+1}$. Between $v_{i-1}d_i^{-1}$ and $v_i$ there will be a vertical bar and a bullet occupying position $\code_i^*(u)+1$, with the elements below displaced to lower rows. Subsequent permutations in the sequence up to the next arrow have their elements displaced around previous bullets as though those rows had been removed. If a number below the bullet in some column is equal to the number immediately to the left of the bar in the same row, then this number will be circled, and this contributes a factor of $y_a-z_b$, where $a$ is the region number (counting from 0 beginning at the left, increasing by $1$ as we pass each bar) and $b$ is the circled number. If $\phi_m(w)\neq 1$, there will be an arrow between $v_m$ and $v_m\phi_m(w)^{-1}$ labeled with $\phi_m(w)$; otherwise, $v_m=v_m\phi_m(w)^{-1}$ will be the final permutation in the sequence. If $v_m\phi_m(w)^{-1}\neq 1$, the diagram contributes a factor of $\sch_{v_m\phi_m(w)^{-1}}(\shiftvby{m}y;z)$, and we write the permutation in bold red type in a column labeled $\sch$. We omit elements $v_i(j)$ of the permutation in the region if $v_i(j)=j$, unless the element is circled.

The significance of the bullets is that they are the positions in the polynomial from which the variable is being pulled out. The relationship between $v_{i-1}d_i^{-1}$ and $v_i$ is that $v_{i-1}d_i^{-1}\downvar{\code_i^*(u)+1}v_i$, so that this is, in a sense, a path of permutations and the bullet occupies the position of $n+1$.

\begin{example} \label{example:figure1}
The permutation $u=[4,3,1,2]$ is dominant, and $\code(u)=(3,2)$. Also, $\code^*(u)=(2,2,1)$. Let $v=[3,1,4,5,2]$ and let $w=[6,3,1,4,2,5]$. Then Theorem \ref{theorem:dualpieri} lets us compute the coefficient $c_{u,v}^w(y;z)$. We have that $\code^*(w)=(2,3,1,1,1)$.
\begin{center}
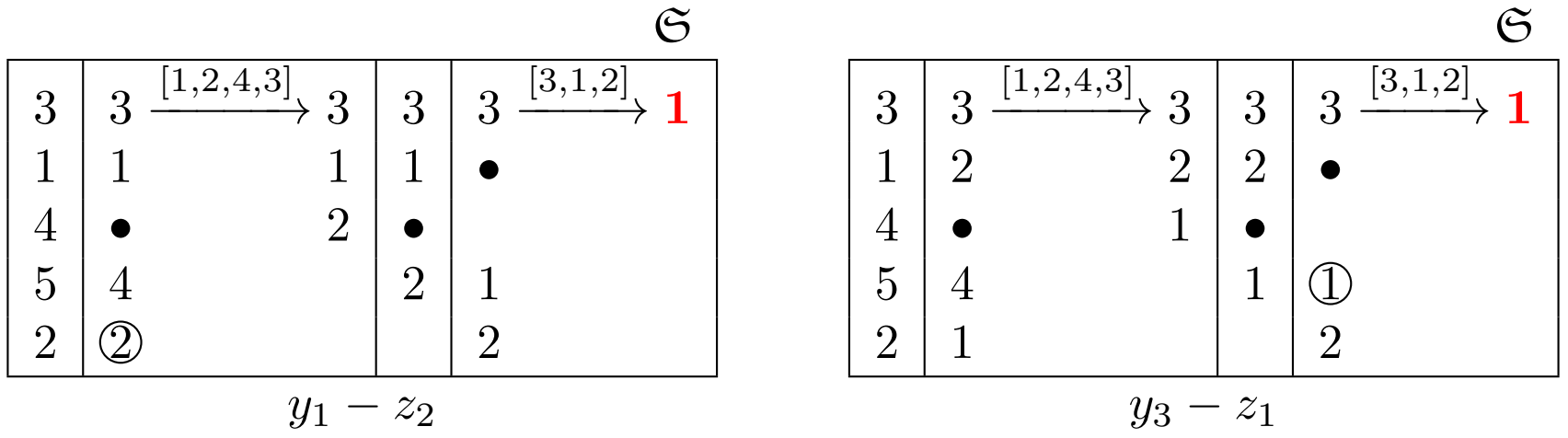

\begin{tabular}{|c|c@{}c@{}c|c|c@{}c@{}c|}
 \multicolumn{8}{r}{$\sch$}\\
\hline
3&3              &\xr{1,2,4,3}&3&3        &3        &\xr{3,1,2}&\rd{1}\\
1&1              &                              &1&1        &\tb&                              &    \\
4&\tb      &                              &2&\tb&         &                              &    \\
5&4              &                              & &2        &1        &                              &           \\
2&\mc{2}  &                              & &         &2        &                              &        \\
\hline
\multicolumn{8}{c}{\makebox[0pt]{$y_1-z_2$}}
\end{tabular}
\hspace{20pt}
\begin{tabular}{|c|c@{}c@{}c|c|c@{}c@{}c|}
\multicolumn{8}{r}{$\sch$}\\
\hline
3&3              &\xr{1,2,4,3}&3&3        &3            &\xr{3,1,2}&\rd{1} \\
1&2              &                              &2&2        &\tb    &                                 & \\
4&\tb      &                              &1&\tb&             &                                 & \\
5&4              &                              & &1        &\mc{1}&                                 &\\
2&1              &                              & &         &2            &                                 &\\
\hline
\multicolumn{8}{c}{\makebox[0pt]{$y_3-z_1$}}
\end{tabular}
\captionof{figure}{Diagrams for the first coefficient in Example~\ref{example:figure1}}
\label{fig:figure1}
\end{center}
Thus $c_{u,v}^w(y;z)=(y_1-z_2)+(y_3-z_1)$. If we take $w=[6,3,1,2,4,5]$ instead, then $\code^*(w)=(2,2,1,1,1)$. Thus, we have the following diagrams.

\begin{center}
\begin{tabular}{|c|c|c|c@{}c@{}c|}
\multicolumn{6}{r}{$\sch$}\\
\hline
3&3              &3               &3        &\xr{3,1,2}&\rd{1} \\
1&1              &1               &\tb&                                 & \\
4&\tb      &                &         &                                 &    \\
5&4              &\tb&         &                                 &           \\
2&\mc{2}  &\mc{2}   &1        &                                 &        \\
 &               &                &2        &                                 &\\
\hline
\multicolumn{6}{c}{\makebox[0pt]{$(y_1-z_2)(y_2-z_2)$}}
\end{tabular}
\hspace{20pt}
\begin{tabular}{|c|c|c|c@{}c@{}c|}
\multicolumn{6}{r}{$\sch$}\\
\hline
3&3              &3               &3                    &\xr{3,1,2}&\rd{1} \\
1&1              &2               &\tb            &                                 & \\
4&\tb      &                &                     &                                 &    \\
5&4              &\tb        &                     &                                 &           \\
2&\mc{2}  &1               &\mc{1}        &                                 &        \\
 &               &                &2                    &                                 &\\
\hline
\multicolumn{6}{c}{\makebox[0pt]{$(y_1-z_2)(y_3-z_1)$}}
\end{tabular}
\hspace{20pt}
\begin{tabular}{|c|c|c|c@{}c@{}c|}
\multicolumn{6}{r}{$\sch$}\\
\hline
3&3              &3               &3                    &\xr{3,1,2}&\rd{1} \\
1&2              &2               &\tb            &                                 & \\
4&\tb      &                &                     &                                 &    \\
5&4              &\tb       &                     &                                 &           \\
2&1              &\mc{1}   &\mc{1}        &                                 &        \\
 &               &                &2                    &                                 &\\
\hline
\multicolumn{6}{c}{\makebox[0pt]{$(y_2-z_1)(y_3-z_1)$}}
\end{tabular}

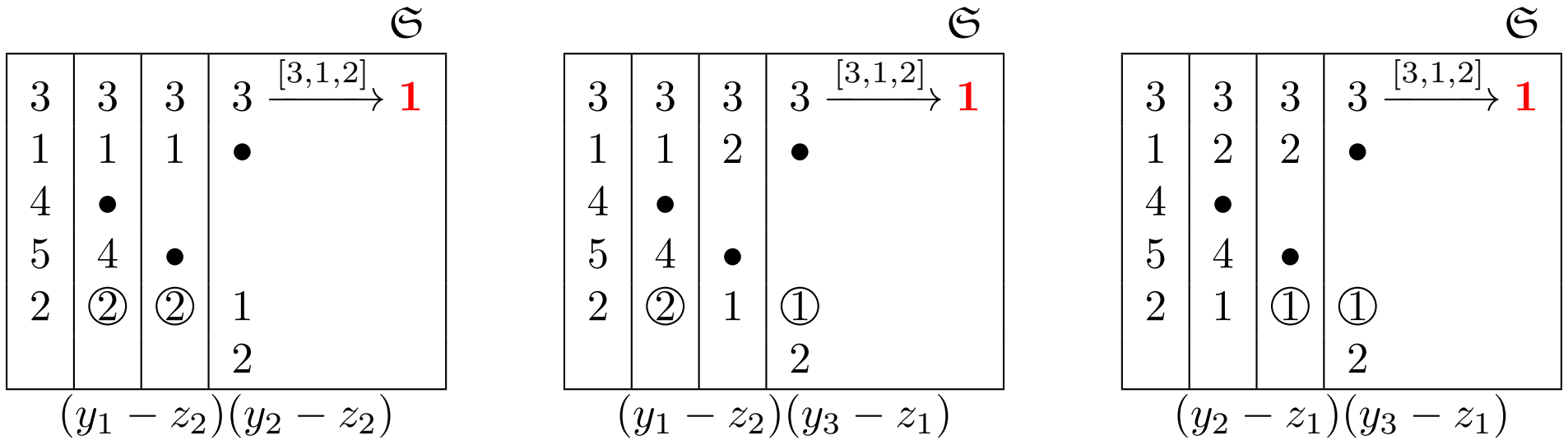
\captionof{figure}{Diagrams for the second coefficient in Example~\ref{example:figure1}}
\label{fig:figure11}
\end{center}

Thus $c_{u,v}^w(y;z)=(y_1-z_2)(y_2-z_2)+(y_1-z_2)(y_3-z_1)+(y_2-z_1)(y_3-z_1)$.

\end{example}

\begin{example} \label{example:figure2}
Let $u=[4,1,2,3]$ and let $v=[1,3,5,2,4]$. We compute $c_{u,v}^w(y;z)$ where $w=[4,1,5,2,3]$. We have
$$\code^*(w)=(1,2,2)$$
$$\code^*(u)=(1,1,1)$$
\begin{center}
\begin{tabular}{|c|c@{}c@{}c|c@{}c@{}c|c|}
\multicolumn{8}{r}{$\sch$}\\
\hline
1&1            &\xr{1,3,2}&1&1        &\xr{1,3,2}&1&\rd{1}\\
3&\tb    &                             &2&\tb&                             &&\rd{\bullet}\\
5&4            &                             &4&3        &                             &&\\
2&\mc{2}&                             &3&2        &                             &&\\
4&3            &                             & &         &                             &&\\
\hline
\multicolumn{8}{c}{\makebox[0pt]{$y_1-z_2$}}
\end{tabular}
\hspace{20pt}
\begin{tabular}{|c|c@{}c@{}c|c@{}c@{}c|c|}
\multicolumn{8}{r}{$\sch$}\\
\hline
1&2            &\xr{1,3,2}&2&2        &\xr{1,3,2}&2&\rd{2}\\
3&\tb    &                             &1&\tb&                             &1&\rd{\bullet}\\
5&4            &                             &4&3        &                             &&\rd{1}\\
2&1            &                             &3&1        &                             &&\\
4&3            &                             & &         &                             &&\\
\hline
\multicolumn{8}{c}{\makebox[0pt]{$\sch_{[2,1]}(y_4;z)$}}
\end{tabular}
\captionof{figure}{Diagrams for Example~\ref{example:figure2}}
\label{fig:figure2}
\end{center}

Thus, the result is
$$c_{u,v}^w(y;z)=(y_1-z_2)+(y_4-z_1)$$
\end{example}

\begin{example} \label{example:figure3}
We compute $c_{[4,2,1,3],[1,3,5,2,4]}^{[5,2,3,1,4]}(y;z)$. We have
$$\code^*(w)=(3,1,1,1)$$
$$\code^*(u)=(2,1,1)$$

\begin{center}
\begin{tabular}{|c@{}c@{}c|c|c|c@{}c@{}c|}
\multicolumn{8}{r}{$\sch$}\\
\hline
1&\xr{1,2,4,3}&1&1            &2        &2        &\xr{2,1}&\rd{1}\\
3&                              &3&3            &\tb&         &                              & \\
5&                              &2&\tb    &         &         &                              &\\
2&                              &5&2            &1        &\tb&                              &\\
4&                              &4&\mc{4}&         &1        &                              &\\
\hline
\multicolumn{8}{c}{\makebox[0pt]{$y_1-z_4$}}
\end{tabular}
\hspace{20pt}
\begin{tabular}{|c@{}c@{}c|c|c|c@{}c@{}c|}
\multicolumn{8}{r}{$\sch$}\\
\hline
1&\xr{1,2,4,3}&1&1            &1            &2        &\xr{2,1}&\rd{1}\\
3&                              &3&3            &\tb    &         &                              &\\
5&                              &2&\tb    &             &         &                              &\\
2&                              &5&4            &3            &\tb&                              &\\
4&                              &4&2            &\mc{2}&1        &                              &\\
\hline
\multicolumn{8}{c}{\makebox[0pt]{$y_2-z_2$}}
\end{tabular}
\begin{tabular}{|c@{}c@{}c|c|c|c@{}c@{}c|}
\multicolumn{8}{r}{$\sch$}\\
\hline
1&\xr{1,2,4,3}&1&1            &2        &2            &\xr{2,1}&\rd{1}\\
3&                              &3&3            &\tb&             &                              &\\
5&                              &2&\tb    &         &             &                              &\\
2&                              &5&4            &3        &\tb    &                              &\\
4&                              &4&2            &1        &\mc{1}&                              &\\
\hline
\multicolumn{8}{c}{\makebox[0pt]{$y_3-z_1$}}
\end{tabular}
\hspace{20pt}
\begin{tabular}{|c@{}c@{}c|c|c|c@{}c@{}c|}
\multicolumn{8}{r}{$\sch$}\\
\hline
1&\xr{1,2,4,3}&1&1            &2            &2        &\xr{2,1}&\rd{1}\\
3&                              &3&4            &\tb    &         &                              &\\
5&                              &2&\tb    &             &         &                              &\\
2&                              &5&2            &1            &\tb&                              &\\
4&                              &4&3            &\mc{3}&1        &                              &\\
\hline
\multicolumn{8}{c}{\makebox[0pt]{$y_2-z_3$}}
\end{tabular}
\begin{tabular}{|c@{}c@{}c|c|c|c@{}c@{}c|}
\multicolumn{8}{r}{$\sch$}\\
\hline
1&\xr{1,2,4,3}&1&1            &3        &3        &\xr{2,1}&\rd{1}\\
3&                              &3&4            &\tb&         &                              &\rd{3}\\
5&                              &2&\tb    &         &         &                              &\rd{2}\\
2&                              &5&2            &1        &\tb&                              &\\
4&                              &4&3            &2        &1        &                              &\\
 &                              & &             &         &2        &                              &\\
\hline
\multicolumn{8}{c}{\makebox[0pt]{$\sch_{[1,3,2]}(y_4,y_5;z)$}}
\end{tabular}

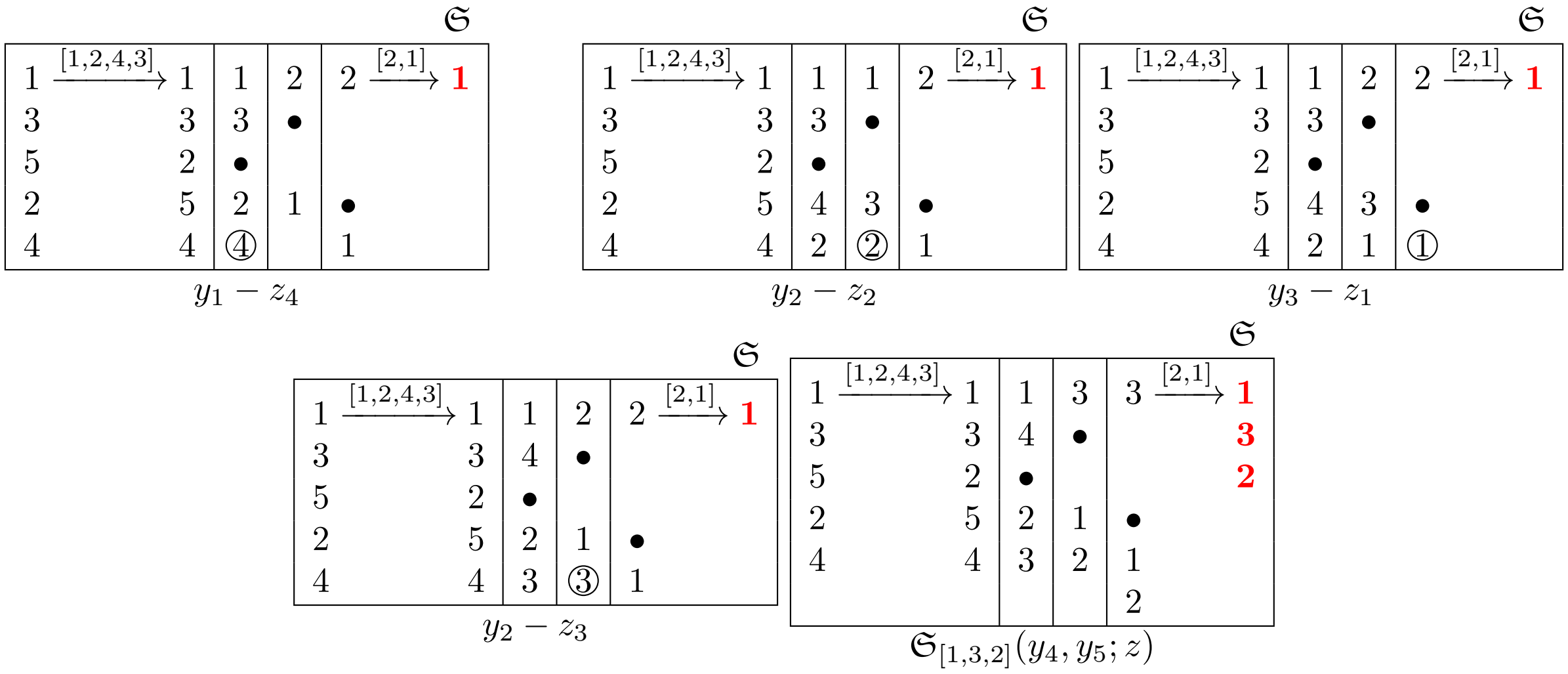
\captionof{figure}{Diagrams for Example~\ref{example:figure3}}
\label{fig:figure3}
\end{center}

Thus, the result is $c_{u,v}^w(y;z)=y_{1} + 2 y_{2} + y_{3} + y_{4} + y_{5} - 2 z_{1} - 2 z_{2} - z_{3} - z_{4}$.
\end{example}

\begin{example} \label{example:figure4}
We give an example where $u$ is not dominant. Set $u=[2, 3, 5, 1, 4]$, set $v=[2, 5, 3, 4, 7, 1, 6]$, and set $w=[2, 6, 3, 4, 7, 1, 5]$. Then
$$\code^*(w)=(5,0,1,1,2)$$
$$\code^*(u)=(3,0,0,1)$$
Then $u$ dominates $w$ but is not dominant. We have the following diagrams.

\begin{center}
\begin{tabular}{|c@{}c@{}c|c|c@{}c@{}c|c|c@{}c@{}c|}
\multicolumn{11}{r}{$\sch$}\\
\hline
2&\xr{1,2,3,5,6,4}&2&2        &\tb    &\xr{2,1}   &2&\tb    &         &\xr{2,3,1}&\rd{1}\\
5&                             &5&5        &3            &                         &3&2            &2&                       &\\
3&                             &3&3        &2            &                         &1&\mc{1}&\tb        &                       &\\
4&                             &1&1        &\mc{1}&                         &5&4            &3        &                       &\\
7&                             &4&\tb&             &                         &4&3            &1        &                       &\\
1&                             &7&6        &5            &                         & &             &         &                       &\\
6&                             &6&4        &\mc{4}&                         & &             &         &                       &\\
\hline
\multicolumn{11}{c}{\makebox[0pt]{$(y_2-z_1)(y_2-z_4)(y_3-z_1)$}}
\end{tabular}
\hspace{2pt}
\begin{tabular}{|c@{}c@{}c|c|c@{}c@{}c|c|c@{}c@{}c|}
\multicolumn{11}{r}{$\sch$}\\
\hline
2&\xr{1,2,3,5,6,4}&2&2        &\tb    &\xr{2,1}    &2&\tb    &         &\xr{2,3,1}&\rd{1}\\
5&                             &5&5        &4            &                         &4&2            &2&                       &\\
3&                             &3&3        &2            &                         &1&\mc{1}&\tb        &                       &\\
4&                             &1&1        &\mc{1}&                         &5&4            &3        &                       &\\
7&                             &4&\tb&             &                         &3&\mc{3}&1        &                       &\\
1&                             &7&6        &5            &                         & &             &         &                       &\\
6&                             &6&4        &3            &                         & &             &         &                       &\\
\hline
\multicolumn{11}{c}{\makebox[0pt]{$(y_2-z_1)(y_3-z_1)(y_3-z_3)$}}
\end{tabular}
\end{center}

\begin{center}
\begin{tabular}{|c@{}c@{}c|c|c@{}c@{}c|c|c@{}c@{}c|}
\multicolumn{11}{r}{$\sch$}\\
\hline
2&\xr{1,2,3,5,6,4}&2&2        &\tb    &\xr{2,1}   &2&\tb    &         &\xr{2,3,1}&\rd{1}\\
5&                             &5&5        &4            &                         &4&3            &3        &                       &\rd{3}\\
3&                             &3&3        &2            &                         &1&\mc{1}&\tb&                       &\rd{2}\\
4&                             &1&1        &\mc{1}&                         &5&4            &2        &                       &\\
7&                             &4&\tb&             &                         &3&2            &1        &                       &\\
1&                             &7&6        &5            &                         & &             &         &                       &\\
6&                             &6&4        &3            &                         & &             &         &                       &\\
\hline
\multicolumn{11}{c}{\makebox[0pt]{$(y_2-z_1)(y_3-z_1)\sch_{[1,3,2]}(y_5,y_6;z)$}}
\end{tabular}
\captionof{figure}{Diagrams for Example~\ref{example:figure4}}
\label{fig:figure4}
\end{center}

Thus
$$c_{u,v}^w(y;z)=(y_2 - z_1) (y_2 - z_4) (y_3 - z_1) + (y_2 - z_1) (y_3 - z_1) (y_3 - z_3) + (y_2 - z_1) (y_3 - z_1)(y_5 + y_6 - z_1 - z_2) $$
\end{example}

\begin{example} \label{example:figure5}
We now compute $c_{[5,4,3,2,1],[4,1,3,5,2]}^{[7,5,3,2,1,4,6]}(y;z)$. We have
$$\code^*(w)=(4,3,2,2,1,1)$$
$$\code^*(u)=(4,3,2,1)$$

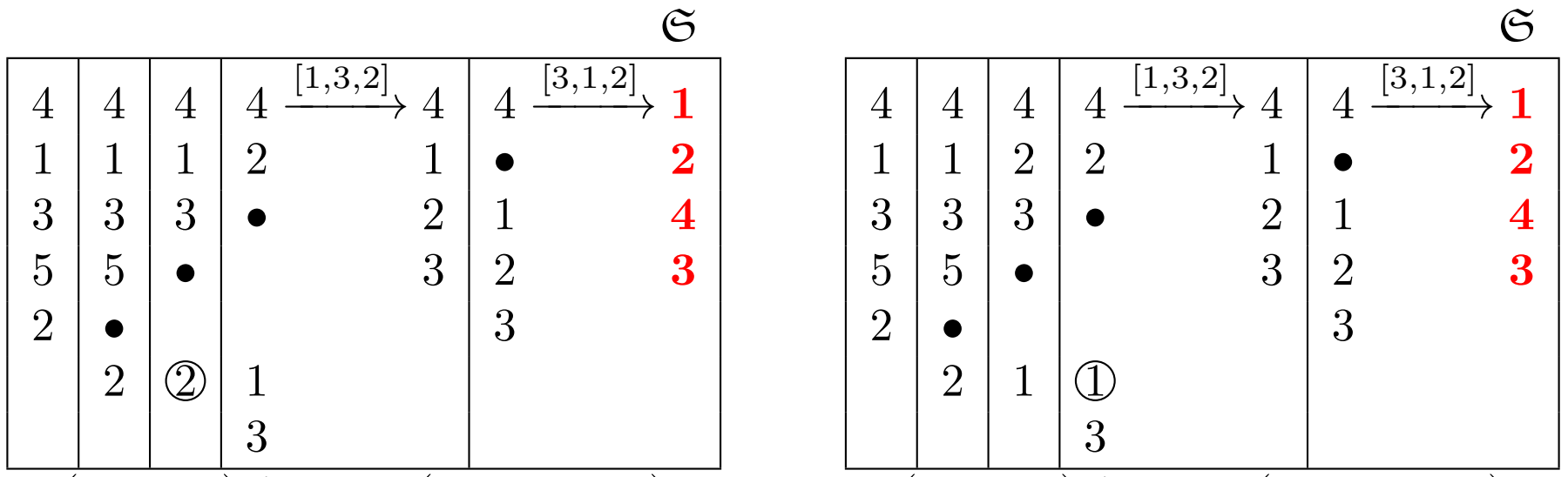
\begin{center}
\begin{tabular}{|c|c|c|c@{}c@{}c|c@{}c@{}c|}
\multicolumn{9}{r}{$\sch$}\\
\hline
4 & 4   & 4     & 4   & \xr{1,3,2} & 4 & 4   & \xr{3,1,2} & \rd{1} \\
1 & 1   & 1     & 2   &            & 1 & \tb &            & \rd{2} \\
3 & 3   & 3     & \tb &            & 2 & 1   &            & \rd{4} \\
5 & 5   & \tb   &     &            & 3 & 2   &            & \rd{3} \\
2 & \tb &       &     &            &   & 3   &            &        \\
  & 2   &\mc{2} & 1   &            &   &     &            &        \\
  &     &       & 3   &            &   &     &            &        \\
\hline
\multicolumn{9}{c}{$(y_2-z_2)\sch_{[1,2,4,3]}(y_5,y_6,y_7;z)$}
\end{tabular}
\hspace{20pt}
\begin{tabular}{|c|c|c|c@{}c@{}c|c@{}c@{}c|}
\multicolumn{9}{r}{$\sch$}\\
\hline
4& 4   & 4   & 4     & \xr{1,3,2} &4 & 4   & \xr{3,1,2} &\rd{1}\\
1& 1   & 2   & 2     &            &1 & \tb &            &\rd{2}\\
3& 3   & 3   & \tb   &            &2 & 1   &            &\rd{4}\\
5& 5   & \tb &       &            &3 & 2   &            &\rd{3}\\
2& \tb &     &       &            &  & 3   &            &\\
 & 2   & 1   & \mc{1}&            &  &     &            &\\
 &     &     & 3     &            &  &     &            &\\
\hline
\multicolumn{9}{c}{\makebox[0pt]{$(y_3-z_1)\sch_{[1,2,4,3]}(y_5,y_6,y_7;z)$}}
\end{tabular}
\captionof{figure}{Diagrams for Example~\ref{example:figure5}}
\label{fig:figure5}
\end{center}

Thus
$$c_{u,v}^w(y;z)=(y_2-z_2)(y_5+y_6+y_7-z_1-z_2-z_3)+(y_3-z_1)(y_5+y_6+y_7-z_1-z_2-z_3)$$
\end{example}

\subsection{Dominant permutations, the code, and manipulation of variables in Schubert polynomials}

Here we gather some results about dominant permutations, the code, the dual code, and the function $\phi_i$ that will be useful throughout this article.

\begin{definition}[Right weak order $\leq_R$] \label{definition:weakorder}
For $u,v\in S_\infty$, we write $u\leq_R v$ if there exists $x\in S_\infty$ such that $v=ux$ and $\ell(v)=\ell(u)+\ell(x)$. For $u\in S_\infty$, let
$$I^*(u)=\{(a,b)\mid a<b\mbox{ and }u^{-1}(a)>u^{-1}(b)\}$$
be the set of pairs of values that appear in decreasing order in the one-line notation of $u$. Then $u\leq_R v$ if and only if $I^*(u)\subseteq I^*(v)$ \cite[Proposition~3.1.3]{combcox}.
\end{definition}

If $\code_i(v)<\code_{i+1}(v)$, then $v(i)<v(i+1)$, so $\ell(vs_i)=\ell(v)+1$. Hence each step in the recursive definition of $\dom(v)$ increases the length by $1$, and $v\leq_R\dom(v)$.

\begin{proposition} \label{proposition:smallestdom}
Let $u\in S_\infty$, and let $\mu$ be a dominant permutation such that $u\leq_R \mu$. Then $\dom(u)\leq_R\mu$. In other words, $\dom(u)$ is the unique minimal dominant permutation above $u$ in right weak order.
\end{proposition}
\begin{proof}
We use induction on $\ell(\dom(u))-\ell(u)$. If $u$ is dominant, then $\dom(u)=u$ and there is nothing to prove. Otherwise, let $i$ be the maximal index such that $\code_i(u)<\code_{i+1}(u)$, so that $u(i)<u(i+1)$ and $\dom(u)=\dom(us_i)$. The indices $j>i+1$ with $u(j)<u(i)$ form a subset of the indices $j>i+1$ with $u(j)<u(i+1)$. The first set has $\code_i(u)$ elements and the second has $\code_{i+1}(u)$ elements, so there exists $j>i+1$ with $u(i)<u(j)<u(i+1)$. Write $a=u(i)$, $b=u(j)$, and $c=u(i+1)$, so that $a<b<c$.

Since $c$ precedes $b$ in the one-line notation of $u$, we have $(b,c)\in I^*(u)\subseteq I^*(\mu)$, so $c$ precedes $b$ in $\mu$. If $a$ preceded $c$ in $\mu$, then $a$, $c$, $b$ would occur in this order in $\mu$ and form a $132$ pattern, contradicting the dominance of $\mu$. Hence $c$ precedes $a$ in $\mu$, so $(a,c)\in I^*(\mu)$. Since
$$I^*(us_i)=I^*(u)\cup \{(a,c)\}$$
we have $I^*(us_i)\subseteq I^*(\mu)$, so $us_i\leq_R\mu$. Since $\ell(\dom(us_i))-\ell(us_i)=\ell(\dom(u))-\ell(u)-1$, the induction hypothesis gives
$$\dom(u) =\dom(us_i) \leq_R \mu$$
\end{proof}

\begin{definition}[$\rho_m$] \label{definition:rhom}
For each $m\geq 0$, define an injective map $\rho_m:S_\infty\to S_\infty$ by inserting $1$ at position $m+1$ and increasing each original entry by $1$. Explicitly, $\rho_m(u)$ is the permutation $u'$ such that
$$u'(i)=\begin{cases}
	u(i)+1&\mbox{ if }i\leq m\\
	u(i-1)+1&\mbox{ if }i>m+1\\
	1&\mbox{ if }i=m+1
\end{cases}$$
Equivalently, $\rho_m(u)=\sigma(u)d[m]$. Its image consists of the permutations $u'$ with $u'(m+1)=1$, and $\rho_m^{-1}$ denotes the inverse map on this image.
\end{definition}

\begin{lemma} \label{lemma:phitrans}
	For $u\in S_\infty$ we have
	$$\phi(u) = \rho_{\code^*_1(u)}^{-1}(u)$$
\end{lemma}
\begin{proof}
For every permutation $v$, we have $\code^*_1(\rho_m(v))=m$, since exactly $m$ entries precede $1$. Inserting $1$ does not change the relative order of the other entries, so $\code^*_{j+1}(\rho_m(v))=\code^*_j(v)$ for all $j\geq 1$. Taking $m=\code^*_1(u)$ and $v=\phi(u)$ gives $\code^*(\rho_m(\phi(u)))=\code^*(u)$. Thus $\rho_m(\phi(u))=u$, proving the result.
\end{proof}

The following is an immediate consequence.

\begin{corollary} \label{corollary:subtractcode}
	Let $u\in S_\infty$ and let $r=\code_1^*(u)+1$. If $j<r$, then
	$$\code_j(\phi(u))=\code_j(u)-1$$
	If $j\geq r$, then
	$$\code_j(\phi(u))=\code_{j+1}(u)$$
\end{corollary}

\begin{lemma} \label{lemma:domdom}
A permutation $u\in S_\infty$ is dominant if and only if for all permutations $w\in S_\infty$ we have that $u$ dominates $w$.
\end{lemma}
\begin{proof}
First, we prove that if $u\in S_\infty$ is dominant, then for any $w\in S_\infty$, we have that $u$ dominates $w$. Let $i$ be the index such that $u(i)=1$, which is $\code_1^*(u)+1$ by Lemma \ref{lemma:phitrans}. Then $\code_i(u)=0$. Since $u$ is dominant, we have that $\code_j(u)=0$ for all $j\geq i$ as well, so that for all $j\geq i$ we have that $u(j)<u(j+1)$. This means that for any permutation $w\in S_\infty$, we have that $u$ $1$-dominates $w$. Now, for any $m\geq 0$ we have that $\phi_m(u)$ is dominant. Thus $\phi_m(u)$ $1$-dominates $\phi_m(w)$ for all $m\geq 0$, regardless of what $w$ is, so that $u$ dominates $w$. Thus, if $u$ is dominant, then $u$ dominates all permutations $w$.

For the converse, we induct on $\ell(u)$, with the identity as the base case. Suppose $u\neq 1$ dominates all permutations. Then $u(j)<u(j+1)$ for all $j>\code_1^*(u)$. In particular, $\code_1^*(u)>0$, and all entries after the position of $1$ are increasing. By definition, $\phi(u)$ dominates all permutations as well. Since $\ell(\phi(u))=\ell(u)-\code_1^*(u)<\ell(u)$, the induction hypothesis shows that $\code(\phi(u))$ is a partition, and we have that
$$\code_j(u)=\begin{cases}
\code_j(\phi(u))+1&\mbox{ if }j\leq\code_1^*(u)\\
0&\mbox{ if }j>\code_1^*(u)
\end{cases}$$
Thus $\code(u)$ is a partition, and the result follows.
\end{proof}

It will be useful to have a formula for $\code^*(\mu)$ in terms of $\code(\mu)$ for $\mu$ dominant.

\begin{lemma} \label{lemma:invdom}
	Suppose $\mu$ is a dominant permutation. Then $\mu^{-1}$ is dominant and $\code^*(\mu)$ is the conjugate partition of $\code(\mu)$, that is to say
	$$\code^*_i(\mu)=\#\{j\mid \code_j(\mu)\geq i\}$$
\end{lemma}
\begin{proof}
	By the discussion in Definition \ref{definition:codestuff}, we have that
	$$\code^*_1(\mu)=\#\{j\mid \code_j(\mu)\geq 1\}$$
	We prove the result by induction on $\ell(\mu)$, the base case of $\ell(\mu)=0$ being clear. We have that
	$$\code^*_i(\phi(\mu))=\code^*_{i+1}(\mu)$$
	for all $i$, and by Corollary \ref{corollary:subtractcode} we have that
	$$\code_i(\phi(\mu))=\code_i(\mu)-1$$
	for all $i$ such that $\code_i(\mu)>0$. By the induction hypothesis, for each $i\geq 1$ we have
	$$\code^*_i(\phi(\mu))=\#\{j\mid \code_j(\mu)-1\geq i\}=\#\{j\mid \code_j(\mu)\geq i+1\}$$
	Thus for $i>1$ we have
	$$\code^*_i(\mu)=\code_{i-1}^*(\phi(\mu))=\#\{j\mid \code_j(\mu)\geq i\}$$
	The result follows by induction.
\end{proof}

\begin{definition}[$\ast$, $\nabla^v$] \label{definition:ydiff}
We wish now to apply divided difference operators to the $y$ variables. We thus define an action $*$ of $S_\infty$ on the polynomial ring permuting the indices of the $y$ variables, and define
$$\nabla^{s_i}(P)=\frac1{y_i-y_{i+1}}(s_i-1) * P$$
Note that we have flipped the sign of the divided difference operator to have the formulas
$$\nabla^{s_i}(\sch_v(x;y))=0$$
if $\ell(s_iv)>\ell(v)$, and
$$\nabla^{s_i}(\sch_v(x;y))=\sch_{s_iv}(x;y)$$
if $\ell(s_iv)<\ell(v)$.
We define skew divided difference operators $\nabla_u^w$ by the analogous identity
$$\nabla^w(PQ)=\sum_u \nabla^u(P)\nabla_u^w(Q)$$
\end{definition}

The basis of our formula is the Pieri formula given in \cite{samuelmolev}, of which we quote a special case.

\begin{proposition}[Special case of the Pieri formula {\cite[Theorem~7.1]{samuelmolev}}] \label{proposition:pieri}
Suppose $k\geq 1$ and $u,w\in S_\infty$. If $u\not\tom{k} w$, then
$$\nabla_u^w\smpr{x_1}{y}{k}=0$$
If $u\tom{k} w$, define
$$Q = \{u(i)\mid i\leq k\mbox{ and }u(i)=w(i)\}$$
then
$$\nabla_u^w\smpr{x_1}{y}{k}=\smpre{x_1}{y_Q}$$
\end{proposition}
\begin{proof}
	This is simply a change of variables from the original theorem.
\end{proof}

For ordinary Schubert polynomials, the next proposition specializes to a formula that isolates a chosen variable $x_i$: one can express $\sch_v(x)$ as a sum of terms of the form $x_i^p\sch_{v'}(x^{(i)})$, leaving Schubert polynomials in the remaining variables \cite[Theorem~5.1]{bsskew}. Proposition \ref{proposition:pullindex} is the double Schubert polynomial version of this result, proved in \cite[Proposition~2.1.5]{samuel2026schubertbialgebra}. We reproduce the proof here for completeness.

\begin{proposition} \label{proposition:pullindex}
Let $v\in S_\infty$ and let $i>0$ be an integer. Then we have
$$\sch_v(x;y)=\sum_{v\downvar{i} v'}\smpre{x_i}{y_{Q_i(v',v)}}\sch_{v'}(x^{(i)};y)$$
\end{proposition}
\begin{proof}
Choose $n\geq i$ with $v\in S_n$, and let $w_0(n)=[n+1,n,\ldots,1]$ be the longest element of $S_{n+1}$. We have
$$\sch_v(x;y)=\nabla^{vw_0(n)}(\sch_{w_0(n)}(x;y))$$
This is equal to
$$\nabla^{vw_0(n)}(\sch_{s_{n+1-i}\cdots s_1w_0(n)}(x^{(i)};y)\smpr{x_i}{y}{n+1-i})$$
and, applying the Leibniz formula, is also equal to
$$\sum_{\substack{v'\in S_\infty\\\ell(v'w_0(n)s_1\cdots s_{n+1-i})=\ell(w_0(n)s_1\cdots s_{n+1-i})-\ell(v')}}\sch_{v'}(x^{(i)};y)\nabla_{v'w_0(n)s_1\cdots s_{n+1-i}}^{vw_0(n)}\smpr{x_i}{y}{n+1-i}$$
By the restricted Pieri formula (Proposition \ref{proposition:pieri}), a summand can be nonzero only if
$$v'w_0(n)s_1\cdots s_{n+1-i}\tom{n+1-i}vw_0(n)$$
This condition is equivalent to
$$v\Tom{i}v'w_0(n)s_1\cdots s_{n+1-i}w_0(n)=v's_ns_{n-1}\cdots s_i$$
We have that $v's_n\cdots s_i$ is exactly $\varphi_{i,n}(v')$ since $v'\in S_n$, so we require that
$$v\Tom{i}\varphi_{i,n}(v')$$
so the sum is over all $v'\in \mathcal{D}_i(v)$.

Applying Proposition \ref{proposition:pieri}, we obtain that the result is equal to
$$\sum_{v'\in \mathcal{D}_i(v)} \smpre{x_i}{y_{A(v',v)}}\sch_{v'}(x^{(i)};y)$$
where $A(v',v)$ is the set of all $vw_0(n)(j)$ such that $1\leq j\leq n+1-i$ and $v'w_0(n)s_1\cdots s_{n+1-i}(j)=vw_0(n)(j)$. For these indices, the equality is equivalent to
$$v'(n+2-s_1\cdots s_{n+1-i}(j))=v(n+2-j)$$
Since $1\leq j\leq n+1-i$, we have
$$s_1\cdots s_{n+1-i}(j)=j+1$$
Hence we need
$$v'(n+2-(j+1))=v'(n+1-j)=v(n+2-j)$$
After substituting $j=n+2-p$, the conditions become $i<p\leq n+1$ and
$$v'(p-1)=v(p)$$
Thus $A(v',v)$ is the set of all $v(p)$ such that $p>i$ and $v'(p-1)=v(p)$, which is exactly $Q_i(v',v)$, and we are done.
\end{proof}

\begin{example}
Suppose $i=2$ and $v=[2,5,1,3,4]$. Then
$$\mathcal{D}_i(v)=\{[2,1],[3,1,2],[4,1,2,3]\}$$
$$Q_i([2,1],[2,5,1,3,4])=\{1,3,4\}$$
$$Q_i([3,1,2],[2,5,1,3,4])=\{1,4\}$$
$$Q_i([4,1,2,3],[2,5,1,3,4])=\{1\}$$
Thus
$$\sch_{[2,5,1,3,4]}(x;y)=\smpre{x_2}{y_1,y_3,y_4}\sch_{[2,1]}(x^{(2)};y)+\smpre{x_2}{y_1,y_4}\sch_{[3,1,2]}(x^{(2)};y)+\smpre{x_2}{y_1}\sch_{[4,1,2,3]}(x^{(2)};y)$$
\end{example}

\begin{example}
Suppose $i=3$ and $v=[2,4,6,5,3,1]$. Then
$$\mathcal{D}_i(v)=\{[2,4,5,3,1],[3,4,5,2,1],[2,5,4,3,1],[3,5,4,2,1]\}$$
$$Q_i([2,4,5,3,1],[2,4,6,5,3,1])=\{1,3,5\}$$
$$Q_i([3,4,5,2,1],[2,4,6,5,3,1])=\{1,5\}$$
$$Q_i([2,5,4,3,1],[2,4,6,5,3,1])=\{1,3\}$$
$$Q_i([3,5,4,2,1],[2,4,6,5,3,1])=\{1\}$$
Thus
\begin{align*}
\sch_{[2,4,6,5,3,1]}(x;y)&=\smpre{x_3}{y_1,y_3,y_5}\sch_{[2,4,5,3,1]}(x^{(3)};y)+\smpre{x_3}{y_1,y_5}\sch_{[3,4,5,2,1]}(x^{(3)};y)\\
                          &+\smpre{x_3}{y_1,y_3}\sch_{[2,5,4,3,1]}(x^{(3)};y)+\smpre{x_3}{y_1}\sch_{[3,5,4,2,1]}(x^{(3)};y)
\end{align*}
\end{example}

\subsection{Relation to pipe dreams}

Proposition \ref{proposition:pullindex} leads to the following positive formula for Schubert polynomials.

\begin{theorem} \label{theorem:schubformula}
	Let $v\in S_\infty$ satisfy $\ell(vs_i)>\ell(v)$ for all $i>n$ for some $n>0$. Then
	$$\sch_v(x;y) = \sum_{\substack{(v_j)_{j=0}^n\\v_0=v,\ v_n=1\\v_{i-1}\downvar{1}v_i\ (1\leq i\leq n)}}\prod_{i=1}^n \smpre{x_i}{y_{Q_1(v_i,v_{i-1})}}.$$
\end{theorem}
\begin{proof}
For $r\geq 0$, write $x^{[r]}=(x_{r+1},x_{r+2},\ldots)$. For $0\leq r<n$, applying Proposition \ref{proposition:pullindex} after replacing each $x_j$ by $x_{j+r}$ gives
$$\sch_{v_r}(x^{[r]};y)=\sum_{v_r\downvar{1}v_{r+1}}\smpre{x_{r+1}}{y_{Q_1(v_{r+1},v_r)}}\sch_{v_{r+1}}(x^{[r+1]};y).$$
Iterating these identities gives the stated sum with an additional factor $\sch_{v_n}(x^{[n]};y)$ in each summand.

For any chain in this expansion, write $q_{i,1},\ldots,q_{i,m_i}$ for the elements of $Q_1(v_i,v_{i-1})$ in decreasing order. Lemma \ref{lemma:factor}, applied at each step and with $\sigma$ applied successively, gives the length-additive factorization
$$v=\left(\prod_{i=1}^n\prod_{j=1}^{m_i}s_{q_{i,j}+i-1}\right)\sigma^n(v_n).$$
If $v_n\neq 1$, it has a right descent $s_j$ for some $j>0$. Then $\sigma^n(v_n)$ has right descent $s_{j+n}$. Since it is a right factor in a length-additive factorization of $v$, $v$ also has right descent $s_{j+n}$, contradicting the hypothesis that $v$ has no right descents beyond $n$. Thus every chain has $v_n=1$, and the residual Schubert polynomial is $1$. This proves the formula.
\end{proof}

\begin{example} \label{example:schubpolyf}
We compute $\sch_{[1,4,3,2]}(x;y)$ with Theorem \ref{theorem:schubformula}.

\begin{center}
\begin{tabular}{|c|c|c|c|}
\hline
1&$\bullet$&       &   \\
4&3      &$\bullet$&    \\
3&2      &\mycircled{2}      &$\bullet$\\
2&1      &\mycircled{1}      &\mycircled{1}\\
\hline
\multicolumn{4}{c}{\makebox[0pt]{$(x_2-y_1)(x_2-y_2)(x_3-y_1)$}}
\end{tabular}
\hspace{64pt}
\begin{tabular}{|c|c|c|c|}
\hline
1&$\bullet$      &        &\\
4&2            &$\bullet$&\\
3&\mycircled{3}&2                 &$\bullet$\\
2&1            &\mycircled{1}      &\mycircled{1}\\
\hline
\multicolumn{4}{c}{\makebox[0pt]{$(x_1-y_3)(x_2-y_1)(x_3-y_1)$}}
\end{tabular}
\hspace{64pt}
\begin{tabular}{|c|c|c|}
\hline
1&$\bullet$      &       \\
4&3            &$\bullet$\\
3&1            &\mycircled{1}\\
2&\mycircled{2}&\mycircled{2}\\
\hline
\multicolumn{3}{c}{\makebox[0pt]{$(x_1-y_2)(x_2-y_1)(x_2-y_2)$}}
\end{tabular}
\end{center}

\begin{center}
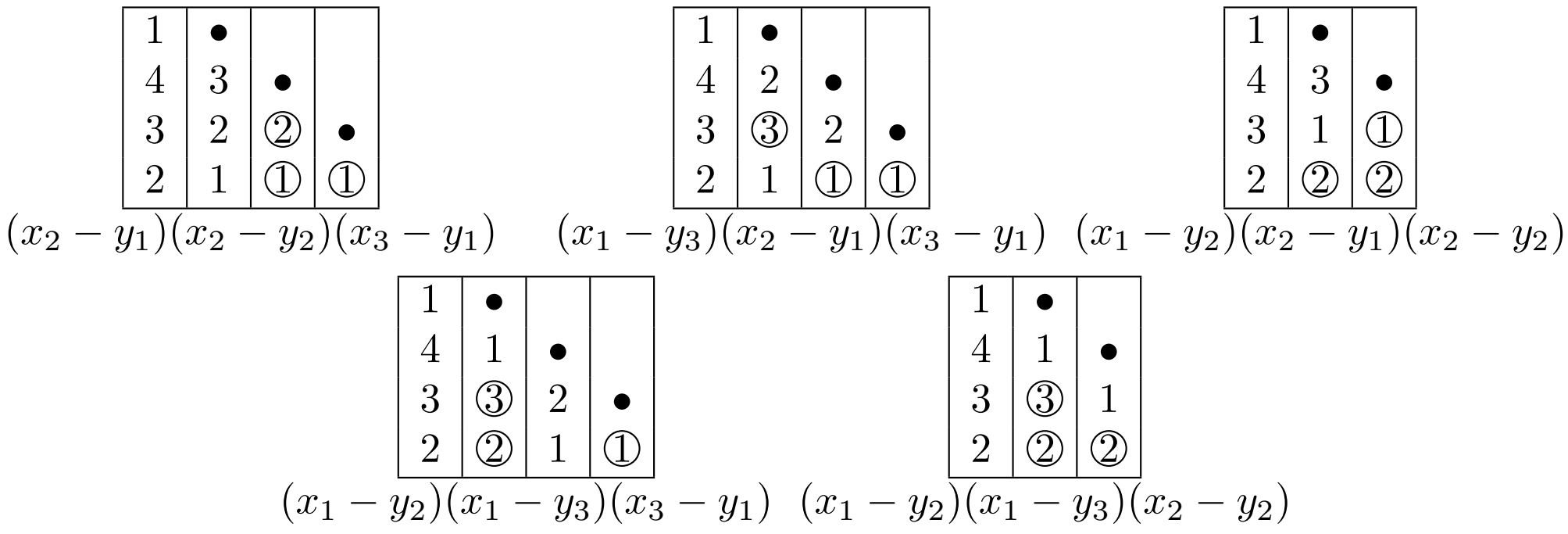

\begin{tabular}{|c|c|c|c|}
\hline
1&$\bullet$      &        &\\
4&1            &$\bullet$&\\
3&\mycircled{3}&2      &$\bullet$\\
2&\mycircled{2}&1      &\mycircled{1}\\
\hline
\multicolumn{4}{c}{\makebox[0pt]{$(x_1-y_2)(x_1-y_3)(x_3-y_1)$}}
\end{tabular}
\hspace{64pt}
\begin{tabular}{|c|c|c|}
\hline
1&$\bullet$      &       \\
4&1            &$\bullet$\\
3&\mycircled{3}&1      \\
2&\mycircled{2}&\mycircled{2}\\
\hline
\multicolumn{3}{c}{\makebox[0pt]{$(x_1-y_2)(x_1-y_3)(x_2-y_2)$}}
\end{tabular}
\captionof{figure}{Diagrams for the double Schubert polynomial in Example~\ref{example:schubpolyf}}
\label{fig:figureschubpolyf}
\end{center}

\end{example}

Theorem \ref{theorem:schubformula} connects the terms in this positive expansion of double Schubert polynomials (products of differences $x_i-y_j$) to reduced words with compatible sequences.

\begin{lemma}\label{lemma:factor}
Let $v$ be a permutation. If $v\downvar{1} v'$, let $a_1,\ldots, a_k$ be the elements of $Q_1(v',v)$ in decreasing order. Then
$$v=s_{a_1}\cdots s_{a_k}\sigma(v')$$
\end{lemma}
\begin{proof}
	Suppose $v\in S_{n}$. We have by definition that there is a sequence of integers $b_1,\ldots,b_p$, all distinct and greater than $1$, such that
	$$vt_{1,b_1}\cdots t_{1, b_p}(1) = n + 1$$
	and
	$$vt_{1,b_1}\cdots t_{1, b_p}(i+1) = v'(i)$$
	for $1\leq i\leq n$. This means that
	$$vt_{1,b_1}\cdots t_{1, b_p} = s_ns_{n-1}\cdots s_1\sigma(v')$$
	This is because if $v''=\sigma(v')$, then
	$$v''(1)=1$$
	and
	$$v''(i) = v'(i-1)+1$$
	The cycle $s_n\cdots s_1$ sends $1\mapsto n+1$ and $i\mapsto i-1$ if $1<i\leq n+1$, hence
	$$vt_{1,b_1}\cdots t_{1, b_p} = s_ns_{n-1}\cdots s_1\sigma(v')$$
	In particular,
	$$v = s_ns_{n-1}\cdots s_1\sigma(v')t_{1,b_p}\cdots t_{1, b_1}$$
	This results in the factorization
	$$v = s_n\cdots s_1t_{1,v'(b_p-1)+1}\cdots t_{1,v'(b_1-1)+1}\sigma(v')$$
	Each forward transposition increases the value at position $1$, so the values $v'(b_1-1),\ldots,v'(b_p-1)$ are strictly increasing. Consequently, multiplying $s_n\cdots s_1$ on the right by these reflections removes the simple reflections $s_{v'(b_p-1)},\ldots, s_{v'(b_1-1)}$. The indices removed are precisely the complement of the elements of $Q_1(v',v)$; hence, the elements of $Q_1(v',v)$ in decreasing order are what remain, as desired.
	The remaining simple reflections are distinct and occur in their original order, so their word is reduced. If the defining $\Tom{1}$ chain has $p$ transpositions, then $\ell(\varphi_{1,n}(v'))=\ell(v')+n=\ell(v)+p$. The $p$ transpositions remove $p$ of the $n$ simple reflections in $s_n\cdots s_1$, leaving $|Q_1(v',v)|=n-p$ reflections. Thus $\ell(v)=\ell(v')+|Q_1(v',v)|$, so the displayed factorization is length-additive.
\end{proof}

\begin{theorem} \label{theorem:sequencecompat}
	Suppose $v=v_0\downvar{1} v_1\cdots \downvar{1} v_{k+1}=1$. For $1\leq i\leq k+1$, let $q_{i,1},\ldots,q_{i,m_i}$ be the elements of $Q_1(v_i,v_{i-1})$ in decreasing order. Let $m=\sum_{i=1}^{k+1}m_i$, and let the sequences $a_1,\ldots,a_m$ and $b_1,\ldots,b_m$ be formed by concatenating, in order of increasing $i$, the sequences of length $m_i$
	$$q_{i,1} + i -1,\ldots, q_{i,m_i} + i - 1$$
	and
	$$(i, i,\ldots, i)$$
	respectively. Then $(a_1,\ldots,a_m)$ is the sequence of indices of the simple reflections in a reduced word for $v$, and $(b_1,\ldots,b_m)$ is a sequence of integers compatible with this reduced word.
\end{theorem}
\begin{proof}
By Lemma \ref{lemma:factor}, each step gives
$$v_{i-1}=\left(\prod_{j=1}^{m_i}s_{q_{i,j}}\right)\sigma(v_i).$$
Applying $\sigma^{i-1}$ to these factorizations and substituting successively yields
$$v=\left(\prod_{i=1}^{k+1}\prod_{j=1}^{m_i}s_{q_{i,j}+i-1}\right)\sigma^{k+1}(v_{k+1})
 =\prod_{i=1}^{k+1}\prod_{j=1}^{m_i}s_{q_{i,j}+i-1},$$
since $v_{k+1}=1$. The length additivity in Lemma \ref{lemma:factor} at each step shows that this concatenated word is reduced.

The sequence $(b_1,\ldots,b_m)$ is weakly increasing, and $b_j\leq a_j$ because each $q_{i,j}\geq1$. Within each block, $b_j$ is constant while the corresponding entries $a_j$ strictly decrease, since the $q_{i,j}$ are listed in decreasing order. These are precisely the compatibility conditions.
\end{proof}

\subsection{Proof of the domination formula}

The idea of the proof is to decompose the divided difference operator and apply it to the double Schubert polynomial one variable at a time.

\begin{lemma} \label{lemma:varshift}
Let $u,w\in S_\infty$, let $f\in\mathbb{Z}[x,y,z]$, and let $i>0$ be an integer. Then
$$\partial_{\sigma_i(u)}^{\sigma_i(w)}(\shiftvby{i}(f))=\shiftvby{i}(\partial_u^w(f))$$
\end{lemma}
\begin{proof}
Since
$$\partial^{s_j}=\frac{1}{x_j-x_{j+1}}(1-s_j)$$
it suffices to show that
$$\sigma_i(s_j)(\shiftvby{i}(x_p))=\shiftvby{i}(s_j(x_p))$$
Writing this out, we have on the left-hand side
$$s_{i+j}(x_{i+p})=\begin{cases}
x_{i+p+1}&\mbox{ if }p=j\\
x_{i+p-1}&\mbox{ if }p=j+1\\
x_{i+p}&\mbox{ otherwise}
\end{cases}$$
which agrees with the right-hand side when worked out.
\end{proof}

\begin{definition}
Let $(a_1,\ldots,a_m)$ be a reduced word. A sequence $(b_1,\ldots,b_m)$ of the same length such that $b_i\in \{1,a_i\}$ is said to be a \emph{subword} for the element $b_1\cdots b_m$. Let $(b_{i_1},\ldots,b_{i_k})$ be all of the nonunit elements of this subword, with $i_1<\cdots<i_k$. Then the subword is said to be reduced if $(b_{i_1},\ldots,b_{i_k})$ is a reduced word.
\end{definition}

The following was proved by Macdonald in \cite{notes}.

\begin{proposition} \label{proposition:divcompute}
Suppose $u,w\in S_\infty$ and let $(a_1,\ldots,a_m)$ be a fixed reduced word for $w$. Then $\partial_u^w$ is the sum over all reduced subwords $(b_1,\ldots,b_m)$ for $u$ of the elements
$$\partial_{b_1}^{a_1}\partial_{b_2}^{a_2}\cdots\partial_{b_m}^{a_m}$$
\end{proposition}
\begin{proof}
This is true if $w=1$. Otherwise, suppose $(a_1,\ldots,a_m)$ is a reduced word for $w$. Let
$$\eta(u,a_m)=\begin{cases}1&\mbox{ if }\ell(ua_m)<\ell(u)\\
0&\mbox{ if }\ell(ua_m)>\ell(u)
\end{cases}$$
Then
$$\partial_u^w=\partial_u^{wa_m}\partial^{a_m}+\eta(u,a_m)\partial_{ua_m}^{wa_m}a_m$$
The set of reduced subwords $(b_1,\ldots,b_m)$ for $u$ is the union of the set of all subwords $(b_1,\ldots,b_{m-1},1)$ that are reduced subwords for $u$ and the set of all subwords $(b_1,\ldots,b_{m-1},a_m)$ such that $(b_1,\ldots,b_{m-1})$ is a reduced subword of $(a_1,\ldots,a_{m-1})$ for $ua_m$, with this second set occurring only if $\eta(u,a_m)=1$. Let
$$D=\partial_{b_1}^{a_1}\cdots\partial_{b_{m-1}}^{a_{m-1}}$$
$D$ occurs in $\partial_{u}^{wa_m}$ where $(b_1,\ldots,b_{m-1})$ is a reduced subword for $u$ by the induction hypothesis, and $\partial^{a_m}=\partial_{1}^{a_m}$, so the terms in the summand $\partial_u^{wa_m}\partial^{a_m}$ correspond bijectively to subwords of the first type. For the second summand, $D$ occurs in $\partial_{ua_m}^{wa_m}$ where $(b_1,\ldots,b_{m-1})$ is a reduced subword for $ua_m$ by the induction hypothesis, and $a_m=\partial_{a_m}^{a_m}$, so the terms in the summand $\partial_{ua_m}^{wa_m}a_m$ correspond bijectively to subwords of the second type. Since this accounts for all of them, the result follows by induction.
\end{proof}

\begin{lemma} \label{lemma:coderesult}
Let $w\in S_\infty$ have $\code^*(w)=(m_1,\ldots,m_q)$. Then
$$w=d[q,m_q]\cdots d[2,m_2]d[m_1]$$
\end{lemma}
\begin{proof}
This is true if $\ell(w)=0$. Otherwise, suppose $\ell(w)>0$ and the result holds for smaller lengths. Assume without loss of generality that $m_q>0$. Then
$$s_qw=d[q+1,m_q-1]d[q-1,m_{q-1}]\cdots d[2,m_2]d[m_1]$$
by the induction hypothesis, since
$$\code^*(s_qw)=(m_1,\ldots,m_{q-1},0,m_q-1)$$
Since $s_qd[q+1,m_q-1]=d[q,m_q]$, the result follows for $s_q(s_qw)=w$ by induction.
\end{proof}

\begin{definition}[The normal form]
The element $d[a,b]$ (for all $a\geq 1$, $b\geq 0$) has a unique reduced word. Consequently, the factorization $w=d[q,m_q]\cdots d[m_1]$ corresponds to an unambiguous reduced word for $w$. We call this reduced word the \emph{normal form} of $w$. This is simply a convenient choice of normal form for our situation. More generally, systematic procedures exist for computing normal forms of elements in Coxeter groups via parabolic subgroups; the construction above is one instance of such methods in the case of symmetric groups \cite{combcox}.
\end{definition}

\begin{lemma} \label{lemma:subword}
Let $(a_1,\ldots,a_m)$ be the normal form for an element $w\in S_\infty$ with $\code^*(w)=(m_1,\ldots,m_q)$, and let $(b_1,\ldots,b_m)$ be a reduced subword for an element $u\in S_\infty$ with $\code^*(u)=(m_1',\ldots,m_p')$, padding either code with zeros as needed. Then $m_1'\leq m_1$. If $(a_j,a_{j+1},\ldots,a_m)$ is the ending portion such that $a_j\cdots a_m=d[m_1]$, then if $m_1'>0$ we have $b_i\neq 1$ for all $j\leq i<j+m_1'$.
\end{lemma}
\begin{proof}
Since $d[m_1]$ and $d[m_1']$ are the minimal-length representatives of the cosets $W_Jw$ and $W_Ju$, respectively, where $W_J$ is the parabolic subgroup generated by $\{s_i\}_{i>1}$, it follows that $d[m_1']\leq d[m_1]$ in Bruhat order. Hence $m_1'\leq m_1$ (see \cite{combcox}). We prove the subword statement by induction on $\ell(u)$.

\emph{Base cases.}
If $m_1'=0$, there is nothing to prove.
If $1$ is the only index $i$ such that $m_i'\neq 0$, then we can see that this is true for combinatorial reasons ($d[m_1']$ has exactly one reduced word, and the only place it fits is at the beginning of $(a_j,\ldots,a_m)$).

\emph{Induction step.}
Suppose $p>1$ is the last index such that $m_p'>0$. Then $\ell(s_pu)<\ell(u)$, so by setting one nonunit element to $1$ in $(b_1,\ldots,b_m)$ we obtain a reduced subword for $s_pu$ (by the exchange property). By the induction hypothesis, the desired portion of $(a_j,\ldots,a_m)$ is occupied by $m_1'$ nonunit elements in this reduced subword for $s_pu$ since $\code_1^*(s_pu)=m_1'$. Replacing the element removed by multiplying on the left by $s_p$, we obtain the result.
\end{proof}

\begin{lemma} \label{lemma:divform}
Let $u, w\in S_\infty$ and suppose $u$ $1$-dominates $w$. Then $\partial_u^w=0$ unless $\code_1^*(w)\geq\code_1^*(u)$, in which case
$$\partial_u^w=\partial_{\sigma(\phi(u))}^{\sigma(\phi(w))}d[\code_1^*(u)]\partial^{d[\code_1^*(u)+1,\code_1^*(w)-\code_1^*(u)]}$$
\end{lemma}
\begin{proof}
Let $(a_1,\ldots,a_m)$ be the normal form of $w$, let $j=m-\code_1^*(w)+1$, and let $(b_1,\ldots,b_m)$ be a reduced subword for $u$. The suffix $(a_j,\ldots,a_m)$ is the reduced word for $d[\code_1^*(w)]$, taken to be empty when $\code_1^*(w)=0$. By Lemma \ref{lemma:subword}, $b_p\neq 1$ for all $j\leq p<j+\code_1^*(u)$ (whereas if $\code_1^*(w)<\code_1^*(u)$, this is not possible and there is no reduced subword for $u$). Since $u$ has no descent at any position $\code_1^*(u)<i\leq\code_1^*(w)$, no further letter of this suffix can be selected: the rightmost such letter would give a descent of $u$. Thus $b_p=1$ for all $p\geq j+\code_1^*(u)$. It follows that
$$\partial_u^w=\partial_{ud[\code_1^*(u)]^{-1}}^{wd[\code_1^*(w)]^{-1}}d[\code_1^*(u)]\partial^{d[\code_1^*(u)+1,\code_1^*(w)-\code_1^*(u)]}$$
since this is true for every subword. Now, we can see by examining the normal forms of these elements that
$$ud[\code_1^*(u)]^{-1}=\sigma(\phi(u))$$
and
$$wd[\code_1^*(w)]^{-1}=\sigma(\phi(w))$$
and the result follows.
\end{proof}

\begin{proof}[Proof of Theorem \ref{theorem:dualpieri}]
We compute this by induction on the maximal index $m$ such that $\code_m^*(u)\neq 0$, with the base case being that $u$ is the identity. If $u$ is the identity, then the result is
$$\partial^w(\sch_v(x;z))|_{x=y}$$
This agrees with the stated formula when we set $m=0$. For the induction step, by Lemma \ref{lemma:divform} we have
$$\partial_u^w=\partial^{\sigma(\phi(w))}_{\sigma(\phi(u))}d[\code^*_1(u)]\partial^{d_1}$$
We apply this to the double Schubert polynomial $\sch_v(x;z)$. Applying $\partial^{d_1}$, we obtain $\sch_{vd_1^{-1}}(x;z)$, provided $\ell(vd_1^{-1})=\ell(v)-\ell(d_1)$, or otherwise $0$. Assuming the term is nonzero, using Proposition \ref{proposition:pullindex}, we have that
$$\sch_{vd_1^{-1}}(x;z)=\sum_{v_1\in \mathcal{D}_{\code_1^*(u)+1}(vd_1^{-1})}\smpre{x_{\code^*_1(u)+1}}{z_{A_1}}\sch_{v_1}(x^{(\code_1^*(u)+1)};z)$$
Now applying $d[\code^*_1(u)]$, this becomes

$$\sum_{v_1\in \mathcal{D}_{\code_1^*(u)+1}(vd_1^{-1})}\smpre{x_1}{z_{A_1}}\sch_{v_1}(\shifthom x;z)$$
We now apply $\partial^{\sigma(\phi(w))}_{\sigma(\phi(u))}$ to this sum. By Lemma \ref{lemma:varshift}, its application to the double Schubert polynomial $\sch_{v_1}(\shifthom x;z)$ is equal to
$$\shifthom\partial^{\phi(w)}_{\phi(u)}(\sch_{v_1}(x;z))$$
We can apply the induction hypothesis to compute this, since $\phi(u)$ dominates $\phi(w)$ and $\code^*(\phi(u))$ is strictly shorter than $\code^*(u)$ (as measured by the position of the rightmost nonzero entry). Therefore, the result follows.
\end{proof}

\subsection{Reversing the separated descents condition}
The article \cite{samuelmolev} uses the properties of dominant permutations and the Pieri formula to compute products of the form $\sch_u(x;y)\sch_v(x;z)$ where there exists a $p$ with $\ell(us_i)>\ell(u)$ if $i<p$ and $\ell(vs_i)>\ell(v)$ if $i>p$. We apply Theorem \ref{theorem:dualpieri} to reverse the conditions in this formula.
\begin{lemma} \label{lemma:shiftup}
Let $u,v,w\in S_\infty$ and let $s_i$ be such that none of $u,v,w$ has $s_i$ as a descent. Then
$$c_{u,v}^w(y;z)=c_{us_i,v}^{ws_i}(y;z)$$
and
$$c_{u,v}^w(y;z)=c_{u,vs_i}^{ws_i}(y;z)$$
\end{lemma}
\begin{proof}
The latter formula was demonstrated in \cite{samuelmolev}. For the former, we apply $\partial_{us_i}^{ws_i}$ to $\sch_v(x;z)$. We have that
$$\partial_{us_i}^{ws_i}=\partial_{us_i}^{w}\partial^{s_i}+\partial_{u}^ws_i$$
We have that $\partial^{s_i}(\sch_v(x;z))=0$, and $s_i(\sch_v(x;z))=\sch_v(x;z)$, so the end result is
$$\partial_u^w(\sch_v(x;z))$$
as desired.
\end{proof}
The following theorem is dual to \cite[Theorem~8.1]{samuelmolev}.
\begin{theorem} \label{theorem:appl}
Let $u,v,w\in S_\infty$. Suppose $u^{-1}\dom(u)=s_{i_1}\cdots s_{i_k}$ is a reduced expression with $\ell(vs_{i_j})>\ell(v)$ for all $1\leq j\leq k$. Then $c_{u,v}^w(y;z)=0$ unless $\ell(wu^{-1}\dom(u))=\ell(w)+\ell(u^{-1}\dom(u))$, in which case
$$c_{u,v}^w(y;z)=c_{\dom(u),v}^{wu^{-1}\dom(u)}(y;z)$$
Hence, this can be computed with Theorem \ref{theorem:dualpieri} and has nonnegative integer coefficients as a polynomial in the terms $y_i-z_j$.

In particular, this applies whenever there exists a $p$ such that $\ell(u)<\ell(us_i)$ for all $i>p$ and $\ell(v)<\ell(vs_i)$ for all $i<p$.
\end{theorem}
\begin{proof}
The argument of \cite[Proposition~4.3]{samuelmolev} applies verbatim.
\end{proof}

Together with \cite[Theorem~3.1]{samuelmolev}, this gives positive combinatorial formulas for $c_{u,v}^w(y;z)$ in the differences $y_i-z_j$ whenever $u$ and $v$ have separated descents in either orientation.
This formula provides the opposite orientation of the generalizations in \cite{samuelmolev}. The formulas in that reference make equivariant positivity manifest, but switching the lower indices is not a symmetry in the general mixed-variable setting. Theorem \ref{theorem:appl} gives positive formulas in the opposite orientation for the generalizations of \cite{grass} (multiplying factorial Schur polynomials with different numbers of $x$ variables) and \cite{kohnert1997multiplication} (multiplying a double Schubert polynomial by a factorial Schur polynomial that has at least as many $x$ variables). Here the factorial Schur polynomial with more $x$ variables has coefficient variables $z$, whereas the earlier formulas used $y$. This also provides a third way to compute a product of two factorial Schur polynomials with the same number of $x$ variables (the first being in \cite{molev1999littlewood}, the second in \cite{samuelmolev}). The resulting formula, however, is not manifestly positive in the simple roots after specializing $z=y$, unlike the formula in \cite{samuelmolev}.
\begin{example} \label{example:dualcomp}
We compute an example dual to \cite[Example~8.1]{samuelmolev}, where $u=[3,1,4,2]$ and $v=[4,1,3,2]$. Here $u$ is not a dominant permutation, but Theorem \ref{theorem:appl} applies. We choose $w=[4,1,5,2,3]$. Then we have
\begin{align*}
\code(\dom(u))&=(2,2)\\
\dom(u)&=[3,4,1,2]\\
wu^{-1}\dom(u)&=[4,5,1,2,3]\\
\code^*(wu^{-1}\dom(u))&=(2,2,2)
\end{align*}

\begin{center}
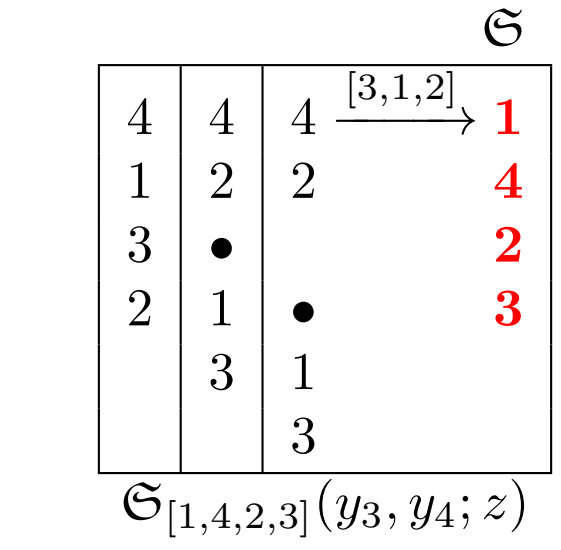

\begin{tabular}{|c|c|c@{}c@{}c|}
\multicolumn{5}{r}{$\sch$}\\
\hline
4&4        &4        &\xr{3,1,2}&\rd{1}\\
1&2        &2        &                       &\rd{4}\\
3&\tb&         &                       &\rd{2}\\
2&1        &\tb&                       &\rd{3}\\
 &3        &1        &                       &\\
 &         &3        &                       &\\
\hline
\multicolumn{5}{c}{\makebox[0pt]{$\sch_{[1,4,2,3]}(y_3,y_4;z)$}}
\end{tabular}
\captionof{figure}{Diagram for Example~\ref{example:dualcomp}}
\label{fig:figuredualcomp}
\end{center}

Thus
$$c_{[3,1,4,2],[4,1,3,2]}^{[4,1,5,2,3]}(y;z)=(y_3-z_2)(y_3-z_3)+(y_3-z_3)(y_4-z_1)+(y_4-z_1)(y_4-z_2)$$
\end{example}

\subsection{Reduction of product coefficients in the case of \texorpdfstring{$1$}{1}-domination}

The induction step in the proof of Theorem \ref{theorem:dualpieri} is, in itself, useful in computing other types of coefficients. The result that is of use is as follows.

\begin{proposition} \label{proposition:oneshift}
	Suppose $u,w\in S_\infty$ are such that $u$ $1$-dominates $w$ and that $v\in S_\infty$. Then $c_{u,v}^w(y;z)=0$ unless $\code_1^*(w)\geq\code_1^*(u)$, so assume this inequality holds. If it is not the case that $$\ell(vd[\code_1^*(u)+1,\code_1^*(w)-\code_1^*(u)]^{-1}) = \ell(v)-(\code_1^*(w) - \code_1^*(u)),$$
then $c_{u,v}^w(y;z)=0$. Otherwise, let $p=\code_1^*(u)$ and let $d=d[\code_1^*(u)+1,\code_1^*(w)-\code_1^*(u)]$. Then
	$$c_{u,v}^w(y;z)=\sum_{vd^{-1}\downvar{p+1}v'}\smpre{y_1}{z_{Q_{p+1}(v',vd^{-1})}}c_{\phi(u),v'}^{\phi(w)}(\shifthom y;z)$$
\end{proposition}
\begin{proof}
	If $c_{u,v}^w(y;z)\neq 0$, then $u\leq w$ in Bruhat order, hence $u^{-1}\leq w^{-1}$ and in particular $u^{-1}(1)\leq w^{-1}(1)$, which is the inequality $\code_1^*(u)\leq\code_1^*(w)$. Applying Lemma \ref{lemma:divform}, we have
	$$\partial_u^w=\partial_{\sigma(\phi(u))}^{\sigma(\phi(w))}s_1s_2\cdots s_p\partial^{d[\code_1^*(u)+1,\code_1^*(w)-\code_1^*(u)]}$$
	The condition on the length of $vd[\code_1^*(u)+1,\code_1^*(w)-\code_1^*(u)]^{-1}$ then arises because otherwise applying this operator yields $0$. Now assume that the condition holds. We have
	$$\sch_{vd^{-1}}(x;z)=\sum_{vd^{-1}\downvar{p+1}v'}\smpre{x_{p+1}}{z_{Q_{p+1}(v',vd^{-1})}}\sch_{v'}(x^{(p+1)};z)$$
	by Proposition \ref{proposition:pullindex}. Applying $s_1\cdots s_p$ yields
	$$s_1\cdots s_p(\sch_{vd^{-1}}(x;z))=\sum_{vd^{-1}\downvar{p+1}v'}\smpre{x_1}{z_{Q_{p+1}(v',vd^{-1})}}\sch_{v'}(\shifthom x;z)$$
	$\partial_{\sigma(\phi(u))}^{\sigma(\phi(w))}$ then acts on $\sch_{v'}(\shifthom x;z)$ as $\partial_{\phi(u)}^{\phi(w)}$ on $\sch_{v'}(x;z)$ by Lemma \ref{lemma:varshift}. The result follows.
\end{proof}

\begin{remark}
	Proposition \ref{proposition:oneshift} is new even for ordinary Schubert polynomials, obtained by specializing $y=z=0$.
\end{remark}

Proposition \ref{proposition:oneshift} enables us to reduce certain product coefficients incrementally and, in some cases, compute the coefficient positively. Indeed, Theorem \ref{theorem:dualpieri} illustrates this method for a wide class of coefficients where full positive computation is achieved. However, other formulas in our framework may still allow positive computation of coefficients even when domination does not apply.

\begin{example} \label{example:pieces}
	Here, it is natural to represent permutations by their codes using the notation $\cperm{\cdots}$. For example, consider computing $c_{\cperm{1,0,3},\cperm{3,0,2,1}}^{\cperm{7,0,1}}(y;z)$. Since $\code_1^*(\cperm{1,0,3}) = 1=\code_1^*(\cperm{7,0,1})$, Proposition \ref{proposition:oneshift} applies, reducing the computation to coefficients of the form $c_{\cperm{0,3},v'}^{\cperm{6,1}}(\shifthom y;z)$. At this stage, however, the prospects are not especially promising, as no special methods are currently available for coefficients of this type.

	Computing $\mathcal{D}_2(\cperm{3,0,2,1})$, we obtain
	$$\{\cperm{3,0,1},\cperm{3,1,1},\cperm{3,2},\cperm{3,2,1}\}$$
	We can compute, using Lemma \ref{lemma:shiftup} on descent position $2$, that
	$$c_{\cperm{0,3},\cperm{3,0,1}}^{\cperm{6,1}}(y;z)=c_{\cperm{0,0,2},\cperm{3,0,1}}^{\cperm{6}}(y;z)=0$$
	since this is a separated descents product, so that \cite[Theorem~3.1]{samuelmolev} applies, and similarly
	$$c_{\cperm{0,3},\cperm{3,1,1}}^{\cperm{6,1}}(y;z)=c_{\cperm{0,0,2},\cperm{3,1,1}}^{\cperm{6}}(y;z)=0$$
	Both $\cperm{3,2}$ and $\cperm{3,2,1}$ are dominant, so these coefficients can be computed with \cite[Theorem~8.1]{samuelmolev}. We have
	$$c_{\cperm{0,3},\cperm{3,2}}^{\cperm{6,1}}(y;z) = y_1 + y_2 - z_1 - z_2$$
	and
	$$c_{\cperm{0,3},\cperm{3,2,1}}^{\cperm{6,1}}(y;z)=(y_1 - z_1) (y_2 - z_1)$$
	We have that $Q_2(\cperm{3,2,1},\cperm{3,0,2,1})$ is empty, and $Q_2(\cperm{3,2},\cperm{3,0,2,1})=\{2\}$.

	The final result is thus
	\begin{align*}
	c_{\cperm{1,0,3},\cperm{3,0,2,1}}^{\cperm{7,0,1}}(y;z)&=(y_1-z_2)c_{\cperm{0,3},\cperm{3,2}}^{\cperm{6,1}}(\shifthom y;z) + c_{\cperm{0,3},\cperm{3,2,1}}^{\cperm{6,1}}(\shifthom y;z)\\
	&=(y_1-z_2)(y_2+y_3-z_1-z_2)+(y_2-z_1)(y_3-z_1)
	\end{align*}

\end{example}

A version of the following result was announced in \cite{mathovschub}. It provides a partial solution to the Schubert--Schur problem: as a consequence of Proposition \ref{proposition:oneshift}, it gives positive combinatorial formulas for Schubert--Schur products whenever the Schur polynomial is ``sufficiently fat.''

\begin{theorem} \label{theorem:schubschurpiece}
	Let $v\in S_\infty$, let $k>0$ be an integer, and let $\lambda$ be a partition with at most $k$ parts. Let $N$ be the maximum value of $0$ and all $\code_i(v)+i-1$ for all $i>k$ such that $\code_i(v)\neq 0$. If $\lambda_i\geq N-k$ for all $i$ with $1\leq i\leq k$, then the coefficient of $\sch_w(x;y)$ in the product
	$$s_{\lambda}(x_1,\ldots,x_k;y)\sch_v(x;z)$$
	admits a positive combinatorial formula in the differences $y_i-z_j$ for all $w\in S_\infty$.
\end{theorem}

\begin{proof}
	Let $u$ be the permutation such that $\sch_u(x;y)=s_\lambda(x_1,\ldots,x_k;y)$. In particular, $\code_i(u)=\lambda_{k+1-i}$ if $1\leq i\leq k$, and $\code_i(u)=0$ if $i>k$. For any permutation $v'$, let $N(v')$ be the maximum value of $0$ and all $\code_i(v')+i-1$ for all $i>k$ such that $\code_i(v')\neq 0$, so that $N=N(v)$. Suppose the polynomial $\sch_w(x;y)$ occurs with a nonzero coefficient in this product.

	We first reduce to the case that $w$ is $k$-Grassmannian. Suppose $w$ has a descent at some $i\neq k$. The product $\sch_u(x;y)\sch_v(x;z)$ is symmetric in $x_i$ and $x_{i+1}$ unless $v$ has a descent at $i$, and in the symmetric case the coefficient of $\sch_w(x;y)$ vanishes because $w$ has a descent at $i$. Hence $v$ has a descent at $i$, and since $u$ does not, Lemma \ref{lemma:shiftup} gives $c_{u,v}^w(y;z)=c_{u,vs_i}^{ws_i}(y;z)$. This replacement does not increase $N(v)$. If $i<k$, then $\code_j(vs_i)=\code_j(v)$ for all $j>k$. If $i>k$, then $\code_i(vs_i)=\code_{i+1}(v)$ and $\code_{i+1}(vs_i)=\code_i(v)-1$, while $\code_j(vs_i)=\code_j(v)$ for $j\notin\{i,i+1\}$. The entry at position $i$ contributes to $N(vs_i)$ only if $\code_{i+1}(v)\neq 0$, and its contribution $\code_{i+1}(v)+i-1$ is less than the contribution $\code_{i+1}(v)+i$ of position $i+1$ to $N(v)$. The contribution of position $i+1$ to $N(vs_i)$ is at most $\code_i(v)+i-1$, which is the contribution of position $i$ to $N(v)$. Since $\ell(ws_i)<\ell(w)$, repeating this process terminates with $w$ $k$-Grassmannian.

	We now prove the result by induction on $N - k$. In the base case $N-k\leq 0$, we have $\code_i(v)=0$ for all $i>k$, so $v$ has no descents after position $k$. Since $u$ has no descents before position $k$, the permutations $u$ and $v$ have separated descents, and the result follows from \cite[Theorem~3.1]{samuelmolev}.

	Suppose $N-k>0$. By virtue of the fact that $c_{u,v}^w(y;z)$ is nonzero, we have that $\code_i(w)\geq \code_i(u)$ for all $1\leq i\leq k$, since Bruhat order on $k$-Grassmannian permutations is containment of partitions. In particular, $\code_i^*(u)=\code_i^*(w)=k$ for all $i$ with $1\leq i\leq N-k$. Applying Proposition \ref{proposition:oneshift}, we have that
	$$c_{u,v}^w(y;z)=\sum_{v'\in \mathcal{D}_{k+1}(v)}\smpre{y_1}{z_{Q_{k+1}(v',v)}}c_{\phi(u),v'}^{\phi(w)}(\shifthom y;z)$$
We have that $\phi(u)$ and $\phi(w)$ are $k$-Grassmannian, and the partition of $\phi(u)$ is $\bar\lambda$ with $\bar\lambda_i=\lambda_i-1$ for all $i$.

We claim that $N(v')\leq N(v)-1$ for all $v'\in\mathcal{D}_{k+1}(v)$. By definition, $\varphi_{k+1,n}(v')=vt_{a_1b_1}\cdots t_{a_pb_p}$, where $a_j\leq k+1<b_j$ and each transposition increases the length by exactly $1$. Right multiplication by such a transposition $t_{ab}$ does not increase $\code_j$ for any $j>k+1$: positions $j>b$ are unaffected, the value at position $b$ decreases, and no position $j$ with $a<j<b$ holds a value strictly between the two values exchanged. Since $\code_i(v')=\code_{i+1}(\varphi_{k+1,n}(v'))$ for all $i\geq k+1$, we obtain $\code_i(v')\leq\code_{i+1}(v)$ for all $i\geq k+1$. Hence, whenever $i>k$ and $\code_i(v')\neq 0$, we have $\code_{i+1}(v)\neq 0$ and
$$\code_i(v')+i-1\leq \code_{i+1}(v)+(i+1)-1-1\leq N(v)-1.$$
Since $N(v)-1\geq k>0$, the claim follows. Therefore
$$\bar\lambda_i=\lambda_i-1\geq N(v)-1-k\geq N(v')-k$$
for all $i$ with $1\leq i\leq k$. Hence, the induction hypothesis applies to $c_{\phi(u),v'}^{\phi(w)}(y;z)$ for each $v'$ occurring in this sum. It follows that all of the terms in the sum are positive; hence, the result follows by induction.
\end{proof}

\begin{remark}
	Theorem \ref{theorem:schubschurpiece} is new even for ordinary Schubert--Schur products.
\end{remark}

\section{Product--splitting equivalence} \label{section:coprod}

\subsection{The equivalence theorem}

Proposition \ref{proposition:pullindex} describes how a double Schubert polynomial expands after splitting its variables. We now express this expansion in terms of double Schubert polynomials in the two sets of variables.

\begin{definition}[Splitting coefficients] \label{definition:coprod}
Let $A$ be a set of positive integers and let $B$ be its complement. We define the $\mathbb{Z}[y,z]$-linear variable-splitting map
$$\operatorname{sp}_A:\mathbb{Z}[x,y,z]\to\mathbb{Z}[x,y,z]\otimes_{\mathbb{Z}[y,z]}\mathbb{Z}[x,y,z]$$
by expressing a polynomial $P(x;y,z)$ in terms of $x_A$ in the first factor and $x_B$ in the second factor. That is, if
$$P(x;y,z)=\sum_i c_i(y,z)P_i(x_A)Q_i(x_B)$$
then set
$$\operatorname{sp}_A(P)=\sum_i c_i(y,z)P_i\otimes Q_i.$$
We define the splitting coefficients $\cpcf{u}{v}{w}{A}{y;z}$ by the expansion
$$\sch_w(x;y)=\sum_{u,v} \cpcf{u}{v}{w}{A}{y;z}\sch_u(x_A;y)\sch_v(x_B;z)$$
When $A$ is finite, the sum is restricted to $u$ with no descents greater than $|A|$, and similarly for $B$ and $v$. These restrictions give Schubert bases in the respective variable sets; coefficients outside these ranges are defined to be zero.

In the ordinary Schubert polynomial case ($y=z=0$), these splitting coefficients are often called coproduct coefficients in the literature. They have a geometric interpretation as the coefficients of the map on cohomology induced by the embedding $Fl(V)\times Fl(W)\hookrightarrow Fl(V\oplus W)$ \cite{coprod}. In equivariant cohomology, an analogous coproduct for back-stable double Schubert polynomials, induced by a direct sum morphism of infinite-dimensional flag varieties, is studied in \cite[\S8]{anderson2025infinite} and \cite[\S9]{anderson2023strong}.
\end{definition}

We will express the splitting coefficients explicitly in terms of triple structure constants. To begin, we give a definition:

\begin{definition}[$\mu_A$] \label{definition:mua}
Let $\mu$ be a dominant permutation, and let $A$ be a set of positive integers (finite or infinite). Indexing $A$ in increasing order and writing it as a function on the corresponding integer interval, we then define $\mu_A$ to be the dominant permutation such that
$$\code_i^*(\mu_A)=\code_{A(i)}^*(\mu)$$
where we define $\code_i^*(\mu_A)=0$ if $i>|A|$.
\end{definition}

\begin{theorem}[Product--splitting equivalence] \label{theorem:doublecoprodpos}
Suppose $w\in S_\infty$ and let $\mu$ be a dominant permutation such that $\ell(w\mu)=\ell(\mu)-\ell(w)$. Let $A$ be a set of positive integers and let $B=\mathbb{N}\setminus A$. Then
$$\sch_w(x;y)=\sum_{\substack{u:\ell(u\mu_A)=\ell(\mu_A)-\ell(u)\\v:\ell(v\mu_B)=\ell(\mu_B)-\ell(v)}}c_{u\mu_A,v\mu_B}^{w\mu}(-y;-z)\sch_u(x_A;y)\sch_v(x_B;z)$$
That is to say, fixing $\mu$, we have that $\cpcf{u}{v}{w}{A}{y;z}=0$ unless $\ell(u\mu_A)=\ell(\mu_A)-\ell(u)$ and $\ell(v\mu_B)=\ell(\mu_B)-\ell(v)$, in which case
$$\cpcf{u}{v}{w}{A}{y;z}=c_{u\mu_A,v\mu_B}^{w\mu}(-y;-z)$$
\end{theorem}

\begin{remark}
	Theorem \ref{theorem:doublecoprodpos} is new even for ordinary Schubert polynomials. Bergeron and Sottile also identify the splitting coefficients of Schubert polynomials with Schubert structure constants \cite[Theorem~1.2.2]{coprod}, but with the structure constants of different products. For $u\in S_n$ and $v\in S_m$, let $u\times v\in S_{n+m}$ denote the permutation that acts as $u$ on $\{1,\ldots,n\}$ and as $v$, shifted by $n$, on $\{n+1,\ldots,n+m\}$. By \cite[Theorem~4.6.1]{coprod}, for each $w$ there are integers $n,m$ and an $n$-element set $\bar{A}\subseteq\{1,\ldots,n+m\}$, obtained from $A$ by truncation and padding, such that $\cpcf{u}{v}{w}{A}{0;0}=0$ unless $u\in S_n$ and $v\in S_m$, in which case
	$$\cpcf{u}{v}{w}{A}{0;0}=c^{(u\times v)\pi}_{\pi,w}(0;0).$$
	Here $\pi\in S_{n+m}$ maps the elements of $\bar{A}$ in increasing order onto $1,\ldots,n$ and the remaining elements in increasing order onto $n+1,\ldots,n+m$. The right-hand side is the coefficient of $\sch_{(u\times v)\pi}$ in the product $\sch_\pi\sch_w$. In this formula the polynomial $\sch_w$ being split is one of the two factors, the other factor $\sch_\pi$ depends only on $\bar{A}$, and $u$ and $v$ appear only in the index of the term of the product. Setting $y=z=0$ in Theorem \ref{theorem:doublecoprodpos} gives instead
	$$\cpcf{u}{v}{w}{A}{0;0}=c^{w\mu}_{u\mu_A,v\mu_B}(0;0),$$
	the coefficient of $\sch_{w\mu}$ in the product $\sch_{u\mu_A}\sch_{v\mu_B}$. The roles of the indices are reversed: $u$ and $v$ index the two factors, and $w$ determines the output index rather than either input index. Comparing the two results yields the identity
	$$c^{(u\times v)\pi}_{\pi,w}(0;0)=c^{w\mu}_{u\mu_A,v\mu_B}(0;0)$$
	for all $u\in S_n$ and $v\in S_m$ satisfying $\ell(u\mu_A)=\ell(\mu_A)-\ell(u)$ and $\ell(v\mu_B)=\ell(\mu_B)-\ell(v)$. Adeyemo and Sottile extended the formula of Bergeron and Sottile to sections of the Billey--Braden pattern map \cite[Theorem~2.5]{adeyemo2014cohomological} and to equivariant cohomology \cite[Theorem~3.3]{adeyemo2015equivariant}. In type $A$ their coefficients have the same form $c^{(u\times v)\pi}_{\pi,w}$, with $\pi$ the minimal length element of a right coset of $S_n\times S_m$ in $S_{n+m}$, so the class being split is again one of the two factors.
\end{remark}

\begin{example}
	We provide a simple example. Set $A=\{1\}$ and $B=\mathbb{N}\setminus A$; only $x_1$ and $x_2$ occur in this example. We have
	\begin{align*}
	\sch_{[1,4,2,3]}(x;y) &= (z_1-y_2)(z_1-y_3)+(z_1-y_3)\sch_{[2,1]}(x_1;y)+(z_1+z_2-y_2-y_3)\sch_{[2,1]}(x_2;z)\\
	&+\sch_{[2,1]}(x_1;y)\sch_{[2,1]}(x_2;z)+\sch_{[3,1,2]}(x_1;y)+\sch_{[3,1,2]}(x_2;z)
	\end{align*}
	For the dominant permutation, we can choose $\mu=\cperm{2,2}$. Then $w\mu = \cperm{1,1}$. We have $\mu_A=\mu_B=\cperm{1,1}$. To find the first coefficient, we compute
	$$c_{\cperm{1,1},\cperm{1,1}}^{\cperm{1,1}}(-y;-z)=(z_1-y_2)(z_1-y_3)$$
	For the second, $u=[2,1]$, and $u\mu_A=\cperm{0,1}$.
	$$c_{\cperm{0,1},\cperm{1,1}}^{\cperm{1,1}}(-y;-z)=z_1-y_3$$
	For the third term, $v=[2,1]$, so
	$$c_{\cperm{1,1},\cperm{0,1}}^{\cperm{1,1}}(-y;-z)=z_1+z_2-y_2-y_3$$
	 For the fourth term, we have
	$$c_{\cperm{0,1},\cperm{0,1}}^{\cperm{1,1}}(-y;-z)=1$$
	For the fifth term, $u=[3,1,2]$, so $u\mu_A$ is the identity, so that the coefficient is $1$, and similarly for the last term where $v=[3,1,2]$.

	The reader might have noticed that the right-hand side has $z$ variables, but the left-hand side does not. In fact, all of the $z$ variables cancel. Nevertheless, the splitting coefficients are polynomials in the differences $z_i-y_j$ with nonnegative integer coefficients, as established below.
\end{example}

The product--splitting equivalence also transfers positivity results to splitting coefficients. We first address the equivariant case $y=z$. By \cite[Theorem~3.3]{adeyemo2015equivariant}, the pullback of an equivariant Schubert class along the embedding $Fl(V)\times Fl(W)\hookrightarrow Fl(V\oplus W)$ expands in the equivariant Schubert basis with equivariant Schubert structure constants as coefficients, which are therefore Graham-positive (see the remark following Theorem \ref{theorem:doublecoprodpos}). In that expansion the equivariant parameters are divided between the two factors, whereas in the corollary below both factors use the full set of parameters $y$. Anderson proved geometrically, by applying Graham's theorem to the image of a direct sum morphism of infinite-dimensional flag varieties, that the closely related equivariant coproduct coefficients of back-stable double Schubert polynomials are Graham-positive \cite[Theorem~8.5]{anderson2025infinite}. See \cite[Theorem~8.7]{anderson2025infinite} for a version with two sets of equivariant parameters, analogous to the mixed-variable case below, and \cite[Theorem~B]{anderson2023strong} for a sharper positivity statement. Those results split the back-stable symmetric part, with a Grassmannian first factor; the following corollary treats the splitting $x=x_A\cup x_B$ for an arbitrary set $A$.

\begin{corollary}[Positivity of double Schubert polynomial splitting coefficients]
Let $w\in S_\infty$, let $A,B\subseteq \mathbb{N}$ be complementary, and write
$$\sch_w(x;y)=\sum_{u,v}\cpcf{u}{v}{w}{A}{y;y}\sch_u(x_A;y)\sch_v(x_B;y)$$
Then the coefficient $\cpcf{u}{v}{w}{A}{y;y}$ is a polynomial in the simple positive roots $y_i-y_{i+1}$ with nonnegative integer coefficients.
\end{corollary}
\begin{proof}
We have that
$$\cpcf{u}{v}{w}{A}{y;y}=c_{u\mu_A,v\mu_B}^{w\mu}(-y;-y)$$
for some $\mu$, since we can always pick some dominant permutation $\mu$ satisfying $\ell(w\mu)=\ell(\mu)-\ell(w)$. By the main result of \cite{graham2001positivity}, $c_{u\mu_A,v\mu_B}^{w\mu}(y;y)$ is a polynomial in the simple negative roots $y_{i+1}-y_i$ with nonnegative integer coefficients. Substitution of $-y$ into this gives the result.
\end{proof}

\begin{corollary}[Positivity of mixed-variable double Schubert polynomial splitting coefficients]
	Let $w\in S_\infty$, let $A,B\subseteq \mathbb{N}$ be complementary, and write
	$$\sch_w(x;y)=\sum_{u,v}\cpcf{u}{v}{w}{A}{y;z}\sch_u(x_A;y)\sch_v(x_B;z)$$
	Then the coefficient $\cpcf{u}{v}{w}{A}{y;z}$ is a polynomial in the terms $z_i-y_j$ with nonnegative integer coefficients.
\end{corollary}
\begin{proof}
	By Theorem \ref{theorem:doublecoprodpos}, each nonzero splitting coefficient has the form $c_{u\mu_A,v\mu_B}^{w\mu}(-y;-z)$ for a suitable dominant permutation $\mu$. Gao and Xiong's positivity theorem \cite[Theorem~1.1]{gao2025grahampositivitytripleschubert} states that $c_{u\mu_A,v\mu_B}^{w\mu}(y;z)$ is a polynomial in the differences $y_i-z_j$ with nonnegative integer coefficients. Substituting $(y,z)\mapsto(-y,-z)$ sends each $y_i-z_j$ to $z_j-y_i$, proving the claim.
\end{proof}

The proof of Theorem \ref{theorem:doublecoprodpos} is an algebraic consequence of the Cauchy formula. The Cauchy formula for double Schubert polynomials is well known \cite{notes}, and related analysis in the context of the nil-Coxeter algebra and the generating function for Schubert polynomials can be found in \cite{fominstan}. The following proposition was proved in greater generality by Kirillov \cite{kirillov2004cauchy} for the universal double Schubert polynomials defined by Fulton \cite{fulton1999universal}. That result specializes immediately to double Schubert polynomials. For completeness, we give a direct elementary proof in the case needed here.

\begin{proposition}[The generalized Cauchy formula for double Schubert polynomials]
Let $w\in S_\infty$ be a permutation and consider an additional infinite set of variables $a=\{a_i\mid i\geq 1\}$. Then we have
$$\sch_w(x;y)=\sum_{\substack{w_1w_2=w\\\ell(w_1)+\ell(w_2)=\ell(w)}}\sch_{w_2}(x;a)\sch_{w_1}(a;y)$$
\end{proposition}
\begin{proof}
The double Schubert polynomials $\sch_u(x;a)$ form a basis of $\mathbb{Z}[x,y,a]$ as a module over $\mathbb{Z}[y,a]$. The expansion of $\sch_w(x;y)$ in this basis is
$$\sch_w(x;y)=\sum_{u\in S_\infty}\partial^u(\sch_w(x;y))|_{x=a}\sch_u(x;a)$$
We have that
$$\partial^u(\sch_w(x;y))|_{x=a}=\begin{cases}
\sch_{wu^{-1}}(a;y)&\mbox{ if }\ell(wu^{-1})=\ell(w)-\ell(u)\\
0&\mbox{ otherwise.}
\end{cases}$$
Thus
$$\sch_w(x;y)=\sum_{\substack{u\in S_\infty\\\ell(wu^{-1})=\ell(w)-\ell(u)}}\sch_{wu^{-1}}(a;y)\sch_u(x;a)$$
This is clearly equivalent to the result.
\end{proof}

We will also need the following factorization property for double Schubert polynomials associated with dominant permutations.

\begin{lemma} \label{lemma:domprod}
Suppose $A$ and $B$ are complementary sets of positive integers, and let $\mu$ be dominant. Then
$$\sch_\mu(x;y)=\sch_{\mu_A}(x;y_A)\sch_{\mu_B}(x;y_B)$$
\end{lemma}
\begin{proof}
Since $A$ and $B$ are complementary, we have that
$$\sch_{\mu}(x;y)=\prod_{a\in A} E_{\code^*_{a}(\mu)}(x;y_a)\prod_{b\in B} E_{\code^*_{b}(\mu)}(x;y_b)$$
By definition of $\mu_A$ and $\mu_B$, and the characterization of double Schubert polynomials corresponding to dominant permutations as products of factorial elementary symmetric polynomials in \cite[Lemma~4.5]{samuelmolev}, the result follows.
\end{proof}

\begin{proof}[Proof of Theorem \ref{theorem:doublecoprodpos}]
By Lemma \ref{lemma:domprod}, we have
$$\sch_\mu(x;y)=\sch_{\mu_A}(x;y_A)\sch_{\mu_B}(x;y_B)$$
Let $a$ and $b$ be additional sets of variables. Combining this identity with the Cauchy formula, we obtain
$$\sch_\mu(x;y)=\sum_{u,v}\sch_{u\mu_A}(x;a)\sch_{v\mu_B}(x;b)\sch_{u^{-1}}(a;y_A)\sch_{v^{-1}}(b;y_B)$$
Here and in the following sums, $u,v$ satisfy $\ell(u\mu_A)=\ell(\mu_A)-\ell(u)$ and $\ell(v\mu_B)=\ell(\mu_B)-\ell(v)$.
Multiplying this out, we obtain
\begin{align*}
\sch_\mu(x;y)&=\sum_w\sum_{u,v}\sch_{w\mu}(x;a)c_{u\mu_A,v\mu_B}^{w\mu}(a;b)\sch_{u^{-1}}(a;y_A)\sch_{v^{-1}}(b;y_B)\\
             &=\sum_w\sum_{u,v}\sch_{w\mu}(x;a)c_{u\mu_A,v\mu_B}^{w\mu}(a;b)\sch_{u}(-y_A;-a)\sch_{v}(-y_B;-b)\\
						 &=\sum_w\sch_{w\mu}(x;a)\sch_{w}(-y;-a)
\end{align*}
In the last line, the Cauchy sum is restricted to $w$ with $\ell(w\mu)=\ell(\mu)-\ell(w)$. For such $w$, comparing the coefficients of $\sch_{w\mu}(x;a)$ in the last two expressions gives
$$\sch_w(-y;-a)=\sum_{u,v}c_{u\mu_A,v\mu_B}^{w\mu}(a;b)\sch_u(-y_A;-a)\sch_v(-y_B;-b)$$
Substituting $x$ for $-y$, $y$ for $-a$, and $z$ for $-b$ gives
$$\sch_w(x;y)=\sum_{u,v}c_{u\mu_A,v\mu_B}^{w\mu}(-y;-z)\sch_u(x_A;y)\sch_v(x_B;z)$$
as desired.
\end{proof}

To justify calling this an ``equivalence'' theorem, we also need to show that product coefficients can be obtained as splitting coefficients.

\begin{definition}[The left dominant approximation $\ldom(u)$, $\ldom(u,v)$] \label{definition:leftdom}
To any permutation $u$ we associate a dominant permutation $\ldom(u)$, the \emph{left dominant approximation of $u$}, as follows. If $u$ is dominant, set $\ldom(u)=u$. Otherwise, if $i$ is the maximal index such that $\code_i^*(u)<\code_{i+1}^*(u)$, set $\ldom(u)=\ldom(s_iu)$. Note that we have
$$\ldom(u)=\dom(u^{-1})^{-1}$$
Then we define $\ldom(u,v)$ to be the dominant permutation such that
$$\code(\ldom(u,v))=\code(\ldom(u))+\code(\ldom(v))$$
\end{definition}

\begin{lemma} \label{lemma:cansplit}
Suppose $u,v\in S_\infty$. Then there exist complementary subsets of positive integers $A$, $B$ such that
$$\ldom(u,v)_A=\ldom(u)$$
and
$$\ldom(u,v)_B=\ldom(v)$$
\end{lemma}
\begin{proof}
By Lemma \ref{lemma:invdom}, $\code^*(\ldom(u,v))$ is the conjugate of the partition $\code(\ldom(u))+\code(\ldom(v))$. The conjugate of a sum of two partitions is the decreasing rearrangement of the concatenation of their conjugates, so $\code^*(\ldom(u,v))$ is the decreasing rearrangement of the concatenation of $\code^*(\ldom(u))$ and $\code^*(\ldom(v))$. Let $A$ be the finite set of positions occupied by the nonzero parts coming from $\code^*(\ldom(u))$ in this rearrangement, with equal parts assigned arbitrarily, and let $B=\mathbb{N}\setminus A$. Then the nonzero parts of $\code^*(\ldom(u,v))$ in positions in $B$ are those of $\code^*(\ldom(v))$, and all remaining positions in $B$ carry $0$. By Definition \ref{definition:mua}, $\ldom(u,v)_A=\ldom(u)$ and $\ldom(u,v)_B=\ldom(v)$.
\end{proof}

\begin{theorem}[Every product coefficient occurs as a splitting coefficient] \label{theorem:allprod}
Suppose $u,v,w\in S_\infty$. Let $A$ and $B$ be subsets of $\mathbb{N}$ as in Lemma \ref{lemma:cansplit}. Then $c_{u,v}^w(y;z)=0$ unless $\ell(w\ldom(u,v)^{-1})=\ell(\ldom(u,v))-\ell(w)$. In that case, we have that
$$c_{u,v}^w(y;z)=\cpcf{u\ldom(u)^{-1}}{v\ldom(v)^{-1}}{w\ldom(u,v)^{-1}}{A}{-y;-z}$$
\end{theorem}
\begin{proof}
Consider the product
$$\sch_{\ldom(u,v)}(x;y)=\sch_{\ldom(u)}(x;y_A)\sch_{\ldom(v)}(x;y_B)$$
Introduce new sets of variables $a,b$, and write
$$\sch_{\ldom(u)}(x;y_A)\sch_{\ldom(v)}(x;y_B)=\sum_{u',v'}\sch_{u'}(x;a)\sch_{v'}(x;b)\sch_{u'\ldom(u)^{-1}}(-y_A;-a)\sch_{v'\ldom(v)^{-1}}(-y_B;-b)$$
In particular, the coefficient of $\sch_{u\ldom(u)^{-1}}(-y_A;-a)\sch_{v\ldom(v)^{-1}}(-y_B;-b)$ is
$$\sch_u(x;a)\sch_v(x;b)=\sum_{w}c_{u,v}^w(a;b)\sch_w(x;a)$$
Considering the Cauchy formula on $\sch_{\ldom(u,v)}(x;y)$, we obtain
$$\sch_{\ldom(u,v)}(x;y)=\sum_{w':\ell(w'\ldom(u,v)^{-1})=\ell(\ldom(u,v))-\ell(w')}\sch_{w'}(x;a)\sch_{w'\ldom(u,v)^{-1}}(-y;-a)$$
By linear independence of double Schubert polynomials, if $c_{u,v}^w(a;b)\neq 0$ then $\sch_w(x; a )$ must occur as a term in this formula. Thus, $\ell(w\ldom(u,v)^{-1})=\ell(\ldom(u,v))-\ell(w)$. By applying Theorem \ref{theorem:doublecoprodpos}, we obtain that
$$c_{u,v}^w(y;z)=c_{(u\ldom(u)^{-1})\ldom(u,v)_A,(v\ldom(v)^{-1})\ldom(u,v)_B}^{(w\ldom(u,v)^{-1})\ldom(u,v)}(y;z)=\cpcf{u\ldom(u)^{-1}}{v\ldom(v)^{-1}}{w\ldom(u,v)^{-1}}{A}{-y;-z}$$
This proves the result.
\end{proof}

The following result illustrates the equivalence between products and splitting coefficients for permutations with a given number of descents.

\begin{theorem}[Descent product--splitting equivalence] \label{theorem:kstep}
Fix $n\geq 1$ and $K\subseteq\{1,\ldots,n-1\}$. If $w\in S_n$ has descents contained in $K$, then there exist a set of descent positions $K'$ with $|K|=|K'|$, a permutation $w'\in S_\infty$ with descents contained in $K'$, and a subset $A\subseteq \mathbb{N}$ such that, for all $u,v\in S_n$ with descents contained in $K$, there exist $u',v'\in S_\infty$ such that
$$c_{u,v}^w(y;z)=\cpcf{u'}{v'}{w'}{A}{-y;-z}$$

Conversely, if $K'$ is a finite set of descent positions, $w'\in S_\infty$ has descents contained in $K'$, and $A\subseteq \mathbb N$, then there exist a set of descent positions $K$ with $|K|\leq |K'|$ and a permutation $w\in S_\infty$ with descents contained in $K$ such that, for all $u',v'\in S_\infty$ with $\cpcf{u'}{v'}{w'}{A}{y;z}\neq 0$, there exist $u,v\in S_\infty$ with descents contained in $K$ such that the above equality holds.
\end{theorem}
\begin{proof}
Let $S_n^K$ denote the set of permutations in $S_n$ with descents contained in $K$, which are the minimal-length representatives of the cosets $\sigma W_J$ for $\sigma\in S_n$, where $W_J$ is the parabolic subgroup generated by $s_i$ for $i\in\{1,\ldots,n-1\}\setminus K$. Let $\mu_1$ be the longest element of $S_n^K$, and let $\mu$ be such that
$$\code_i(\mu)=2\code_i(\mu_1)$$
for all $i$. The permutations $\mu_1$ and $\mu$ are dominant, with descent set $K$. Every inversion of $t\in S_n^K$ lies between two different consecutive blocks determined by $K$, and both $\mu_1$ and $\mu$ invert every such pair of positions. By the inversion-set characterization of weak order applied to inverses, we have
$$\ell(t\mu_1^{-1})=\ell(\mu_1)-\ell(t),\qquad \ell(t\mu^{-1})=\ell(\mu)-\ell(t).$$
The conjugate partition $\code^*(\mu)$ repeats each part of $\code^*(\mu_1)$ twice, so there exist complementary sets $A,B$ such that $\mu_A=\mu_B=\mu_1$. By Theorem \ref{theorem:doublecoprodpos}, for any $u,v\in S_n^K$, we have that
$$c_{u,v}^w(y;z)=\cpcf{u\mu_1^{-1}}{v\mu_1^{-1}}{w\mu^{-1}}{A}{-y;-z}$$
Let $K'$ be the descent set of $\mu^{-1}$. Since a dominant permutation and its inverse have the same number of descents, $|K'|=|K|$. We may select $w'=w\mu^{-1}$, $u'=u\mu_1^{-1}$, and $v'=v\mu_1^{-1}$. Since $\mu^{-1}=w^{-1}(w\mu^{-1})$ with $\ell(\mu^{-1})=\ell(w^{-1})+\ell(w\mu^{-1})$, every descent of $w\mu^{-1}$ is a descent of $\mu^{-1}$. Thus $w'$ has descents contained in $K'$, and the displayed identity holds for all the prescribed $u,v$, including when the coefficient is zero.

Conversely, suppose $w'$ has descents contained in $K'$. Let $\mu=\ldom(w')^{-1}=\dom((w')^{-1})$, and let $K$ be the descent set of $\mu$. Since $\mu$ is obtained from $(w')^{-1}$ by successive right multiplications by simple transpositions, each of which increases the length by $1$, we have $\ell(w'\mu)=\ell(\mu)-\ell(w')$. By Proposition \ref{proposition:mudescents}, $\desc(\mu^{-1})=\desc(w')$; since $\mu$ is dominant, it has the same number of descents as its inverse. Therefore $|K|=|\desc(w')|\leq |K'|$. Then for any $A$ and $B=\mathbb{N}\setminus A$, $\mu_A$ and $\mu_B$ have descents contained in $K$. For any $u',v'$ with $\cpcf{u'}{v'}{w'}{A}{y;z}\neq 0$, Theorem \ref{theorem:doublecoprodpos} gives
$$\cpcf{u'}{v'}{w'}{A}{-y;-z}=c_{u'\mu_A,v'\mu_B}^{w'\mu}(y;z)$$
and its length conditions imply that all three indices in the product coefficient have descents contained in $K$. Set $w=w'\mu$, $u=u'\mu_A$, and $v=v'\mu_B$, and the result follows.
\end{proof}

Thus, a combinatorial formula for products in $k$-step flag varieties is equivalent to a formula for splitting coefficients of Schubert polynomials indexed by permutations with $k$ descents. In fact, known results imply combinatorial positivity of splitting coefficients of double Schubert polynomials with $z=y$ indexed by permutations with at most two descents, and of ordinary Schubert polynomials with at most three descents.

\begin{corollary}
The puzzle rule in \cite{buch2015mutations} computes $\cpcf{u}{v}{w}{A}{y;y}$ Graham-positively for all $u$, $v$, $w$, and $A$ such that $w$ has at most two descents.
\end{corollary}

\begin{corollary}
The puzzle rule proved in \cite{threestep} computes $\cpcf{u}{v}{w}{A}{0;0}$ positively for all $u$, $v$, $w$, and $A$ such that $w$ has at most three descents.
\end{corollary}

\subsection{A historical missed opportunity: separated descents for ordinary Schubert polynomials}

Let $k\geq 1$ be an integer, and consider the variable-splitting expansion
$$\sch_w(x)=\sum_{u,v}d_w^{u,v}(k)\sch_u(x_1,\ldots,x_k)\sch_v(x_{k+1},x_{k+2},\ldots)$$
The coefficients $d_w^{u,v}(k)$ were computed in 1991 by Macdonald in \cite[(4.19)]{notes} and shown to be nonnegative. This variable-splitting expansion is often called the Schubert polynomial coproduct in the literature. The condition of separated descents was introduced in 2019 in \cite{kzj}, and positive puzzle formulas for this case appear in \cite{kzjsep}. Recall that permutations $u,v$ are said to have separated descents if there exists a $p>0$ such that $\ell(us_i)>\ell(u)$ for all $i<p$ and $\ell(vs_i)>\ell(v)$ for all $i>p$. We now interpret this condition in terms of Theorem \ref{theorem:allprod}. We begin with the following proposition.

\begin{proposition} \label{proposition:mudescents}
Let $u\in S_\infty$ and suppose $j>0$ is an integer. Then $\ell(us_j)<\ell(u)$ if and only if $\ell(\ldom(u)s_j)<\ell(\ldom(u))$.
\end{proposition}
\begin{proof}
We prove this by induction on $\ell(\ldom(u))-\ell(u)$, the result being clear if $u=\ldom(u)$. If $u$ is not dominant, let $i$ be the maximal index such that $\code_i^*(u)<\code_{i+1}^*(u)$. We know by definition that $\ldom(u)=\ldom(s_iu)$. If $\ell(us_j)<\ell(u)$, then $\ell(s_ius_j)<\ell(s_iu)$. By the induction hypothesis, $\ell(\ldom(u)s_j)<\ell(\ldom(u))$.

Conversely, suppose $\ell(us_j)>\ell(u)$ and assume the induction hypothesis. We claim that $\ell(s_ius_j)>\ell(s_iu)$. If instead we assume that $\ell(s_ius_j)<\ell(s_iu)$, then in fact $s_iu=us_j$ by the exchange property. This means that $u^{-1}(i)=j$, $u^{-1}(i+1)=j+1$. However, this implies that $\code_i^*(u)=\code_{i+1}^*(u)$, contradicting the choice of $i$. Thus we can conclude that $\ell(s_ius_j)>\ell(s_iu)$, and hence by the induction hypothesis $\ell(\ldom(u)s_j)>\ell(\ldom(u))$, as desired. This completes the induction and proves the proposition.
\end{proof}

If $u$ and $v$ have separated descents, therefore, so do $\ldom(u)$ and $\ldom(v)$. We also have the following.

\begin{lemma} \label{lemma:separateddom}
Let $\mu$ be a dominant permutation, let $A=[k]$ for some $k>0$, let $B=\mathbb N\setminus [k]$, and let $p=\code^*_k(\mu)$. Then $\ell(\mu_As_i)>\ell(\mu_A)$ for all $i<p$ and $\ell(\mu_Bs_i)>\ell(\mu_B)$ for all $i>p$, or in other words $\mu_A$ and $\mu_B$ have separated descents.
\end{lemma}
\begin{proof}
If $p=0$, then $\mu_B=1$, and the result follows immediately. Otherwise, observe that the nonzero entries of $\code^*(\mu_A)$ are the right descent positions of $\mu_A$, and similarly for $\mu_B$. Since $\code^*(\mu)$ is a partition, $\mu_A$ has no descents $i$ such that $i<p$ because $\code^*_j(\mu_A)\geq p$ for $1\leq j\leq k$, and $\mu_B$ has no descents $i$ such that $i>p$ because $\code_j^*(\mu_B)\leq p$ for all $j$. This establishes the claim.
\end{proof}

Conversely, suppose $\ell(us_i)>\ell(u)$ for all $i<p$ and $\ell(vs_i)>\ell(v)$ for all $i>p$. Then the sets $A$ and $B$ of Lemma \ref{lemma:cansplit} can be chosen to be intervals (one initial and finite, and the other complementary and infinite). After specializing $y=z=0$, the coefficients needed to compute $c_{u,v}^w(0;0)$ using Theorem \ref{theorem:allprod} are exactly those in Macdonald's variable-splitting expansion. Thus, a formula for separated descents could already have been obtained in 1991, 28 years before its appearance in \cite{kzj}.

We now establish positivity of these splitting coefficients in the equivariant and mixed-variable settings.

\begin{theorem}[Separated splitting theorem]
Let $k>0$ and let $A=[k]$ be the initial interval. Then the coefficients
\[
\cpcf{u}{v}{w}{A}{y;z}
\]
appearing in the expansion
\[
\sch_w(x;y) \;=\; \sum_{u,v} \cpcf{u}{v}{w}{A}{y;z}\,
   \sch_u(x_{1},\ldots,x_k;y)\,
   \sch_v(x_{k+1},x_{k+2},\ldots;z),
\]
have the following properties:
\begin{enumerate}
    \item They are polynomials with nonnegative integer coefficients in the differences $z_i-y_j$.
    \item After the specialization $z=y$, they are \emph{Graham-positive}, that is, polynomials in the simple roots $y_i-y_{i+1}$ with nonnegative integer coefficients.
\end{enumerate}
Moreover, nonnegativity in the differences $z_i-y_j$ also holds when $A=\mathbb{N}\setminus[k]$.
\end{theorem}

\begin{proof}

Assume first that $A=[k]$ and set $B=\mathbb{N}\setminus [k]$. We have that $\cpcf{u}{v}{w}{A}{y;z}=c_{u\mu_A,v\mu_B}^{w\mu}(-y;-z)$ for some dominant permutation $\mu$. Since $\mu_A=u^{-1}(u\mu_A)$ with $\ell(\mu_A)=\ell(u^{-1})+\ell(u\mu_A)$, every descent of $u\mu_A$ is a descent of $\mu_A$, and similarly every descent of $v\mu_B$ is a descent of $\mu_B$. By Lemma \ref{lemma:separateddom}, there exists a $p$ such that $\ell(u\mu_As_i)>\ell(u\mu_A)$ for all $i<p$ and $\ell(v\mu_Bs_i)>\ell(v\mu_B)$ for all $i>p$. Hence, \cite[Theorem~3.1]{samuelmolev} allows us to compute this coefficient nonnegatively as a polynomial in $z_i-y_j$ that is manifestly Graham-positive after the specialization $z=y$. For $A=\mathbb{N}\setminus [k]$, the opposite descent conditions hold, and Theorem \ref{theorem:appl} applies to provide positivity as a polynomial in $z_i-y_j$. When $z=y$, the coefficients for $A=\mathbb{N}\setminus[k]$ coincide with those for $A=[k]$ with the roles of $u$ and $v$ exchanged. Graham positivity is manifest, without cancellation, in the formula of \cite[Theorem~3.1]{samuelmolev}, but not in the formula of Theorem \ref{theorem:appl}.
\end{proof}

\subsection{Reduction formulas for product coefficients by varying the dominant permutation}

The splitting coefficient in Theorem \ref{theorem:doublecoprodpos} does not depend on which admissible dominant permutation is used to express it as a product coefficient. Suppose we wish to compute the product coefficient $c_{u,v}^w(y;z)$. Let $A$ and $B$ be as in Lemma \ref{lemma:cansplit}, and write
\begin{align*}
	\widehat{w}&=w\ldom(u,v)^{-1}\\
	\widehat{u}&=u\ldom(u)^{-1}\\
	\widehat{v}&=v\ldom(v)^{-1}
\end{align*}
We may think of $\widehat{u}$, $\widehat{v}$, and $\widehat{w}$ as the ``base indices'' of this splitting coefficient. By Theorem \ref{theorem:allprod}, $c_{u,v}^w(y;z)=0$ unless $\ell(\widehat{w})=\ell(\ldom(u,v))-\ell(w)$, in which case
$$c_{u,v}^w(y;z)=\cpcf{\widehat{u}}{\widehat{v}}{\widehat{w}}{A}{-y;-z}$$
The dominant permutation $\ldom(u,v)$ used to express this splitting coefficient as a product coefficient may, however, be larger than necessary. Let
$$\mu=\dom(\widehat{w}^{-1})$$
Since $\widehat{w}^{-1}\leq_R\mu$, we have $\ell(\widehat{w}\mu)=\ell(\mu)-\ell(\widehat{w})$. Applying Theorem \ref{theorem:doublecoprodpos} with this $\mu$, we find that $c_{u,v}^w(y;z)=0$ unless, in addition, $\ell(\widehat{u}\mu_A)=\ell(\mu_A)-\ell(\widehat{u})$ and $\ell(\widehat{v}\mu_B)=\ell(\mu_B)-\ell(\widehat{v})$, in which case
$$c_{u,v}^w(y;z)=c_{\widehat{u}\mu_A,\widehat{v}\mu_B}^{\widehat{w}\mu}(y;z)$$
This transformation never increases the length of the upper index. Indeed, the condition $\ell(\widehat{w})=\ell(\ldom(u,v))-\ell(w)$ states that $\widehat{w}^{-1}\leq_R\ldom(u,v)$, so $\mu\leq_R\ldom(u,v)$ by Proposition \ref{proposition:smallestdom}, and therefore
$$\ell(\widehat{w}\mu)=\ell(\mu)-\ell(\widehat{w})\leq\ell(\ldom(u,v))-\ell(\widehat{w})=\ell(w)$$
At worst, $\mu=\ldom(u,v)$, and the indices remain unchanged. The following example demonstrates that this method can convert a coefficient that we do not immediately know how to compute into one that we do know how to compute.
\begin{example} \label{example:reduction}
Let
\begin{align*}
	u&=\cperm{1,0,2}\\
	v&=\cperm{0,2,1}\\
	w&=\cperm{4,0,1}
\end{align*}
Then
$$c_{u,v}^w(y;z)=y_1+y_5-z_1-z_2$$
None of the known formulas allows a direct positive computation of this coefficient. The reduction method, however, transforms it into $c_{\cperm{0,0,2},\cperm{0,1,1}}^{\cperm{0,1,2}}(y;z)$, which the reader may verify is equal to the original coefficient. Since
$$\sch_{\cperm{0,1,1}}(x;z)=s_{(1,1)}(x_1,x_2,x_3;z)$$
this case can be computed using the Pieri formula \cite[Theorem~7.1]{samuelmolev}. Alternatively, one may apply the factorial Schur product formula \cite[Theorem~3.1]{molev1999littlewood}, as the expression is a product of factorial Schur polynomials.
\begin{align*}
&s_{(2)}(x_1,x_2,x_3;y)\,s_{(1,1)}(x_1,x_2,x_3;z)\\
&\quad=\bigl((y_1 - z_1) (y_5 - z_2) + (y_2 - z_1) (y_1 - z_1) + (y_2 - z_2) (y_5 - z_2)\bigr)s_{(2)}(x_1,x_2,x_3;y) \\
&\qquad+(y_1 + y_2 - z_1 - z_2)s_{(3)}(x_1,x_2,x_3;y)+(y_1 + y_5 - z_1 - z_2)s_{(2,1)}(x_1,x_2,x_3;y)\\
&\qquad+s_{(3,1)}(x_1,x_2,x_3;y)+s_{(2,1,1)}(x_1,x_2,x_3;y)
\end{align*}
whereas the original product we were trying to compute, $\sch_{\cperm{1,0,2}}(x;y)\sch_{\cperm{0,2,1}}(x;z)$, has $25$ nonzero double Schubert polynomial terms.
\end{example}

The \texttt{schubmult} Python package can display product and splitting coefficients in positive form. For the most general algorithms applicable to arbitrary polynomials, obtaining such a presentation tends to be computationally intractable---even in small cases---when approached directly through integer programming without problem-specific optimizations. By contrast, when the reduction formulas introduced in this paper are applied, together with other positive formulas established in prior work, the package frequently produces positive presentations for cases of reasonable size. Notably, such optimizations were essential to render positivity testing of coefficients feasible in the case $n=6$.

\section*{Acknowledgments}

We thank Kellen Meyers for assistance with large-scale computations used to test our conjectures.

\bibliographystyle{acm}
\bibliography{formschub}

@misc{threestep,
  title={Schubert puzzles and integrability {I}: invariant trilinear forms},
  author={Knutson, Allen and Zinn-Justin, Paul},
  howpublished={arXiv:1706.10019 [math.CO]},
  year={2017}
}

@misc{kzjsep,
  title={Schubert puzzles and integrability {III}: separated descents},
  author={Knutson, Allen and Zinn-Justin, Paul},
  howpublished={arXiv:2306.13855 [math.CO]},
  year={2023}
}

@misc{anderson2023strong,
	title={Strong equivariant positivity for homogeneous varieties and back-stable coproduct coefficients},
	author={Anderson, David},
	howpublished={arXiv:2302.12765 [math.AG]},
	year={2023},
	url={https://arxiv.org/abs/2302.12765}
}

@article{anderson2025infinite,
	title={Infinite flags and {S}chubert polynomials},
	author={Anderson, David},
	journal={Forum Math. Sigma},
	volume={13},
	pages={Paper No. e36, 27 pp.},
	year={2025},
	doi={10.1017/fms.2024.99}
}

@article{adeyemo2014cohomological,
	title={Cohomological consequences of the pattern map},
	author={Adeyemo, Praise and Sottile, Frank},
	journal={Illinois J. Math.},
	volume={58},
	number={4},
	pages={997--1008},
	year={2014}
}

@misc{adeyemo2015equivariant,
	title={Equivariant cohomology theories and the pattern map},
	author={Adeyemo, Praise and Sottile, Frank},
	howpublished={arXiv:1506.04411 [math.AG]},
	year={2015},
	url={https://arxiv.org/abs/1506.04411}
}

@article{laisant1888numeration,
  title={Sur la num{\'e}ration factorielle, application aux permutations},
  author={Laisant, C-A},
  journal={Bull. Soc. Math. France},
  volume={16},
  pages={176--183},
  year={1888}
}

@article{fan2025bumpless,
  title={Bumpless pipe dreams meet puzzles},
  author={Fan, Neil JY and Guo, Peter L and Xiong, Rui},
  journal={Adv. Math.},
  volume={463},
  pages={110113},
  year={2025},
  publisher={Elsevier}
}

@misc{schubmult,
	author = {Matthew J. Samuel},
	title = {{\texttt{schubmult}: a {P}ython package for multiplying {S}chubert polynomials}},
	howpublished = {Available from the {P}ython {P}ackage {I}ndex, \url{https://pypi.org/project/schubmult/}, with sources at \url{https://github.com/matthematics/schubmult}},
	year = {2026}
}

@misc{samuel2026schubertbialgebra,
	author = {Samuel, Matthew J.},
	title = {The {S}chubert bialgebra},
	howpublished = {arXiv:2610.22582 [math.CO]},
	year = {2026},
	url = {https://arxiv.org/abs/2610.22582}
}

@book{humphreys1992reflection,
  title={Reflection groups and {C}oxeter groups},
  author={Humphreys, James E},
  year={1992},
  publisher={Cambridge University Press}
}

@misc{kzj,
  title={Schubert calculus and quiver varieties},
  author={Allen Knutson and Paul Zinn-Justin},
  howpublished={\url{https://pi.math.cornell.edu/~allenk/boston2019.pdf}},
	note="Accessed 2023-06-10",
  year={2019}
}

@article{bsskew, 
	author = "Nantel Bergeron and Frank Sottile",
	journal = "Trans. Amer. Math. Soc.",
	year = "2002",
	volume = "354",
	number = "2",
	title = "Skew {S}chubert functions and the {P}ieri formula for flag manifolds", 
	pages = "651--673"
}

@phdthesis{samuelleibniz, 
	author = "Matthew J. Samuel",
	school = "Rutgers, The State University of New Jersey",
	year = "2014",
	title = "The {L}eibniz formula for divided difference operators associated to {K}ac-{M}oody root systems"
}

@BOOK{combcox,
  TITLE = {Combinatorics of {C}oxeter {G}roups},
  SUBTITLE = {Graduate Texts in Mathematics 231},
  AUTHOR = {Anders Bj{\"o}rner and Francesco Brenti},
  YEAR = {2005}, 
  PUBLISHER = {Springer},
	PLACE = "New York"
}

@article{coprod,
author = "Nantel Bergeron and Frank Sottile",
title = "Schubert polynomials, the {B}ruhat order, and the geometry of flag manifolds",
journal = "Duke Math. J.",
volume = "95",
year = "1998", 
number = "2", 
pages = "373--423"
}

@article{fominstan,
author = "S. Fomin and R. P. Stanley",
title = "Schubert {P}olynomials and the {N}il{C}oxeter {A}lgebra",
journal = "Adv. Math.",
volume = "103",
year = "1994", 
number = "2", 
pages = "196--207"
}

@article{fulton1999universal,
  title={Universal {S}chubert polynomials},
  author={Fulton, William},
  journal={Duke Math. J.},
  volume={96},
  number={3},
  pages={575--594},
  year={1999}
}

@article{kirillov2004cauchy,
  title={Cauchy identities for universal {S}chubert polynomials},
  author={Kirillov, Anatolii Nikolaevich},
  journal={J. Math. Sci. (N.Y.)},
  volume={121},
  pages={2360--2370},
  year={2004},
  publisher={Springer}
}

@article{samuelmolev,
    author =       "Samuel, Matthew J.",
    title =        "A {M}olev-{S}agan type formula for double {S}chubert polynomials",
    journal =      "J. Pure Appl. Algebra",
    year =         "2024",
		volume = "228",
		number = "7",
		pages = "107636"
}

@book{notes,
    author    = "Macdonald, I. G.",
    title     = "Notes on {S}chubert {P}olynomials",
    year      = "1991",
    address = "D{\'e}p. de math{\'e}matiques et d{'}informatique",
    publisher   = "Universit{\'e} du Qu{\'e}bec {\`a} Montr{\'e}al"
}

@article{sottile,
		author = "Sottile, Frank",
		title = "Pieri{'}s rule for flag manifolds and {S}chubert polynomials",
		year = "1996",
		journal = "Ann. Inst. Fourier (Grenoble)",
		volume = "46",
		number = "1",
		pages = "89--110"
}

@article{grass,
		author = "Kevin Purbhoo and Frank Sottile",
		title = "A {L}ittlewood-{R}ichardson rule for {G}rassmannian permutations",
		year = "2009",
		journal = "Proc. Amer. Math. Soc.",
		volume = "137",
		pages = "1875--1882"
}

@article{kohnert1997multiplication,
  title={Multiplication of a {S}chubert polynomial by a {S}chur polynomial},
  author={Kohnert, Axel},
  journal={Ann. Comb.},
  volume={1},
  number={1},
  pages={367--375},
  year={1997},
  publisher={Springer}
}

@article{reif,
title = {On the diagram of 132-avoiding permutations},
journal = {European J. Combin.},
volume = {24},
number = {6},
pages = {759-776},
year = {2003},
author = {Astrid Reifegerste},
}

@article{graham2001positivity,
  title={Positivity in equivariant {S}chubert calculus},
  author={Graham, William},
  journal={Duke Math. J.},
  volume={109},
  number={3},
  pages={599--614},
  year={2001},
  publisher={Duke University Press}
}

@article{buch2015mutations,
  title={Mutations of puzzles and equivariant cohomology of two-step flag varieties},
  author={Buch, Anders Skovsted},
  journal={Ann. of Math.},
  volume={182},
  number={1},
  pages={173--220},
  year={2015},
  publisher={JSTOR}
}

@article{molev1999littlewood,
  title={A {L}ittlewood-{R}ichardson rule for factorial {S}chur functions},
  author={Molev, Alexander and Sagan, Bruce},
  journal={Trans. Amer. Math. Soc.},
  volume={351},
  number={11},
  pages={4429--4443},
  year={1999}
}

@MISC {mathovschub,
    title = {``{S}ufficiently fat'' {S}chur polynomial times {S}chubert polynomial is {S}chubert-positive},
    author = {Samuel, Matthew J.},
    howpublished = {{M}ath{O}verflow},
    note = {\url{https://mathoverflow.net/q/453974} (version: 2023-09-04)},
    year = {2023},
    url = {https://mathoverflow.net/q/453974}
}

@article{gao2025grahampositivitytripleschubert,
	title={Positivity in triple {S}chubert calculus},
	author={Gao, Yibo and Xiong, Rui},
	journal={S{\'e}m. Lothar. Combin.},
	volume={95B},
	note={Proceedings of FPSAC 2026, Art.~37, 11 pp.},
	year={2026},
	url={https://radon.mat.univie.ac.at/~slc/wpapers/FPSAC2026/37.html}
}

\end{document}